\documentclass[a4paper,12pt]{article}
\usepackage{amsmath,amssymb}
\usepackage{eucal}
\usepackage[OT1,T1]{fontenc}
\advance\topmargin by-\headheight
\advance\topmargin by-\headsep
\advance\jot by2pt

\newcommand\rk{\mathop{\mathrm{rk}}\nolimits}
\renewcommand\deg{\mathop{\mathrm{deg}}\nolimits}
\newcommand\im{\mathop{\mathrm{im}}\nolimits}
\newcommand\coker{\mathop{\mathrm{coker}}\nolimits}
\newcommand\Hom{\mathop{\mathrm{Hom}}\nolimits}
\newcommand\Ext{\mathop{\mathrm{Ext}}\nolimits}
\newcommand\End{\mathop{\mathrm{End}}\nolimits}
\newcommand\Aut{\mathop{\mathrm{Aut}}\nolimits}
\newcommand\Iso{\mathop{\mathrm{Iso}}}
\newcommand\Coh{\mathop{\mathrm{Coh}}\nolimits}
\newcommand\Pone{\mathop{\mathbf P^1}}
\newcommand\Fq{{\mathord{\mathbf F_q}}}
\newcommand\tor{{\mathord{\mathrm{tor}}}}
\newcommand\lf{{\mathord{\mathrm{lf}}}}
\newcommand\ad{\mathop{\mathrm{ad}}}
\newcommand\catmodf[1]{#1\mbox{-${\underline{\mathrm{modf}}}$}}
\newcommand\Uqsltwohat{{U_q(\widehat{\mathfrak{sl}_2})}}
\newcommand\longlongrightarrow{\relbar\joinrel\relbar\joinrel\relbar\joinrel
  \relbar\joinrel\relbar\joinrel\relbar\joinrel\relbar\joinrel\rightarrow}
\renewcommand\labelitemi{$\ \bullet$}
\renewcommand\theenumi{(\roman{enumi})}
\renewcommand\labelenumi{\theenumi}

\renewcommand\enumerate{\list\labelenumi
  {\setlength\leftmargin{0pt}
  \usecounter{enumi}\def\makelabel##1{\kern\labelsep{##1}}}}
\renewcommand\itemize{\list\labelitemi
  {\setlength\leftmargin{0pt}\setlength\labelwidth{0pt}
  \def\makelabel##1{\kern\labelsep{##1}}}}
\newtheorem{theorem}{Theorem}
\newtheorem{lemma}[theorem]{Lemma}
\newtheorem{proposition}[theorem]{Proposition}
\newtheorem{corollary}[theorem]{Corollary}
\newenvironment{proof}{\trivlist
  \item[\hskip \labelsep{\itshape Proof.}]\upshape}{\nobreak\noindent$\square$\endtrivlist}
\newenvironment{other}[1]{\refstepcounter{theorem}\trivlist
  \item[\hskip \labelsep{\itshape #1~\arabic{theorem}.}]\upshape}{\endtrivlist\bigbreak}

\begin{document}
\title{The Hall algebra of the category of coherent sheaves on the projective line}
\author{Pierre Baumann and Christian Kassel}
\date{}
\maketitle
\begin{abstract}
To an abelian category $\mathcal A$ of homological dimension 1 satisfying certain
finiteness conditions, one can associate an algebra, called the Hall algebra.
Kapranov studied this algebra when $\mathcal A$ is the category of
coherent sheaves over a smooth projective curve defined over a finite field,
and observed analogies with quantum affine algebras. We recover here in an
elementary way his results in the case when the curve is the projective~line.
\end{abstract}
\begingroup
\small
\begin{center}
\bfseries R\'esum\'e
\end{center}
\quotation
A toute cat\'egorie ab\'elienne $\mathcal A$ de dimension homologique~$1$ v\'erifiant
certaines conditions de finitude, on peut associer une alg\`ebre appel\'ee l'alg\`ebre de Hall.
Kapranov a \'etudi\'e cette alg\`ebre lorsque $\mathcal A$ est la cat\'egorie des faisceaux
coh\'erents sur une courbe projective lisse d\'efinie sur un corps fini et a observ\'e des
analogies entre l'alg\`ebre de Hall et les alg\`ebres affines quantiques. Nous
red\'emontrons de mani\`ere \'el\'ementaire ses r\'esultats dans le cas o\`u la courbe est
la droite projective.
\endquotation
\endgroup

\section*{Introduction}
The combinatorial lattice structure of objects in an abelian category $\mathcal A$ of homological
dimension~$1$ satisfying certain finiteness conditions may be encoded in an algebraic structure,
called the Hall algebra of $\mathcal A$. Hall's original results, as described in Chapters~II and III
of Macdonald's book \cite{Macdonald}, concern the category of modules of finite length over a discrete
valuation ring with finite residue field.

A decade ago Ringel studied the Hall algebra of the category of finite-dimensional
representations over a finite field $\Fq$ of a quiver whose underlying graph $\Gamma$
is a Dynkin diagram of type A, D, or E. He showed that a suitable modification of this Hall
algebra yields an algebra isomorphic to the ``positive part'' of Drinfeld's and Jimbo's
quantized enveloping algebra associated to~$\Gamma$ and specialized at the value $q$ of the parameter
(see \cite{Green2} or Section~2 of~\cite{Ringel3} for an introduction to Ringel's results).

More recently, Kapranov investigated the case when $\mathcal A$ is the category of
coherent sheaves over a smooth projective curve $X$ defined over~$\Fq$. In a remarkable
paper \cite{Kapranov}, he used unramified automorphic forms, Eisenstein series and
$L$-functions to translate the geometrical properties of $X$ into the algebraic structure of
Ringel's modification of the Hall algebra of~$\mathcal A$. This allowed him to observe some striking
similarities between such a Hall algebra and Drinfeld's loop realization of the quantum affine
algebras~\cite{Drinfeld2}. In the case when $X$ is the projective line~$\Pone(\Fq)$,
Kapranov deduced from his general constructions an isomorphism between a
certain subalgebra of the Hall algebra and a certain ``positive part'' of the
(untwisted) quantum affine algebra~$\Uqsltwohat$.

The main objective of this article is to recover Kapranov's isomorphism for the projective
line in a more elementary way. Avoiding any use of adelic theory or of automorphic forms, we
compute directly the structure constants of the Hall algebra of the category of coherent
sheaves over~$\Pone(\Fq)$. This approach moves us away from the analogy that motivated Kapranov,
but hopefully makes his results more accessible and concrete. We also observe that Kapranov's
isomorphism yields a natural definition for the vectors of the Poincar\'e-Birkhoff-Witt basis
of~$\Uqsltwohat$ that Beck, Chari, and Pressley introduced in~\cite{Beck-Chari-Pressley}.

The paper is organized as follows. In Section~\ref{se:HallAlg} we give the definitions of
the Hall algebra and of its Ringel variant associated to an abelian category satisfying
adequate conditions. In Section~\ref{se:CohSheavesP1} we recall basic facts on
the category $\mathcal A$ of coherent sheaves on the projective line~$\Pone(k)$
over an arbitrary field $k$ and we analyse carefully the extensions between certain
``elementary'' objects. This leads in Section~\ref{ss:HallNumCohShP1} to
Theorem~\ref{th:HallNumCohShP1}, which provides many structure constants of the Hall algebra of $\mathcal A$
when~$k=\Fq$. Now every coherent sheaf can be written as the direct sum of its torsion subsheaf
and of a locally free subsheaf. The existence of such decompositions gives rise to a factorization
of the Hall algebra as a semidirect product of two subalgebras, denoted below by $B_1$ and
$H(\mathcal A_\tor)$, and related to locally free coherent sheaves and torsion sheaves, respectively.
By an averaging process which takes into account all closed points of $\Pone(\Fq)$, we define
in Sections~\ref{ss:HallAlgA_tor} and \ref{ss:SubalgKapra} a subalgebra $B_0$ of $H(\mathcal A_\tor)$.
In the final Section~\ref{se:Uqsl2hat}, we recall the definition of the quantum affine
algebra $\Uqsltwohat$ and we relate it to the subalgebra of the Hall algebra generated by $B_0$ and $B_1$.

Our interest in this subject grew out of a seminar held in Strasbourg during the\break
winter 1996--97 and aimed at understanding Kapranov's paper~\cite{Kapranov}.
We are grateful to Henri~Carayol, Florence~Lecomte, Louise~Nyssen, Georges~Papadopoulo,
and Marc~Rosso for their enlightening lectures. The first author also acknowledges the
financial support of the French Minist\`ere de l'Education Nationale, de la Recherche
et de la Technologie and of the~CNRS.

\section{Hall algebras}\label{se:HallAlg}
\subsection{Initial data}\label{ss:InitData}
Let $k$ be a field. Recall that an abelian category $\mathcal A$ is said to be $k$-linear if
the homomorphism groups in $\mathcal A$ are endowed with the structure of $k$-vector
spaces, the composition of morphisms being $k$-bilinear operations. In the sequel, we
will consider abelian $k$-linear categories $\mathcal A$ satisfying the following
finiteness conditions (H1)--(H4):
\begin{description}
\item[(H1)]The isomorphism classes of objects in $\mathcal A$ form a set $\Iso(\mathcal A)$.
\item[(H2)]For all objects $V$, $W$ in $\mathcal A$, the $k$-vector space
$\Hom_{\mathcal A}(V,W)$ is finite-dimensional.
\item[(H3)]For all objects $V$, $W$ in $\mathcal A$, the $k$-vector space
$\Ext^1_{\mathcal A}(V,W)$, defined as a set of equivalence classes of
short exact sequences of the form
\begin{equation*}
0\to W\to U\to V\to0,
\end{equation*}
is finite-dimensional.
\item[(H4)]$\mathcal A$ can be imbedded as a full subcategory in an abelian $k$-linear
category $\mathcal B$ with enough injectives (or projectives), $\mathcal A$ is closed
under extensions in $\mathcal B$, and $\Ext^2_\mathcal B(V,W)=0$ for all objects $V$,
$W$ in $\mathcal A$.
\end{description}
\noindent The isomorphism class of an object $V$ in $\mathcal A$ will be denoted by
$[V]\in\Iso(\mathcal A)$, and the isomorphism class of the zero object
by $0$. It will be convenient to choose a preferred object $M(\alpha)$ in
each isomorphism class $\alpha\in\Iso(\mathcal A)$. Condition (H2) implies that
$\mathcal A$ satisfies the Krull-Schmitt property: each object $V$ in $\mathcal A$ can
be written as a direct sum $W_1\oplus\cdots\oplus W_\ell$ of indecomposable objects,
the isomorphism classes of the objects $W_i$ and their multiplicities in the decomposition
being uniquely determined. Finally, Condition~(H4) ensures that short exact sequences
in $\mathcal A$ give rise to 6-term exact sequences of $k$-vector spaces involving the
bifunctors $\Hom_{\mathcal A}(-,-)$ and $\Ext^1_{\mathcal A}(-,-)$.

Among the categories $\mathcal A$ that we will consider, certain enjoy an additional finiteness
condition (H5), namely:
\begin{description}
\item[(H5)]Each object in $\mathcal A$ has a finite filtration with simple quotients
(Jordan-H\"older series).
\end{description}

\subsection{The Grothendieck group and the Euler form}\label{ss:GrotEuler}
The Grothendieck group $K(\mathcal A)$ of the category $\mathcal A$ is, by definition,
the abelian group presented by the generators $d(\alpha)$, where
$\alpha\in\Iso(\mathcal A)$, together with the relations
$d(\beta)=d(\alpha)+d(\gamma)$ whenever there is a short exact sequence
\begin{equation*}
0\longrightarrow M(\gamma)\longrightarrow M(\beta)\longrightarrow M(\alpha)\longrightarrow0.
\end{equation*}
If $V$ is an object in $\mathcal A$, we will write $d(V)$ instead of $d([V])$
to denote the image of its class in the Grothendieck group. If $\mathcal A$ satisfies
Condition~(H5), then $K(\mathcal A)$ is the free abelian group on the symbols $d(\alpha)$,
for all isomorphism classes $\alpha\in\Iso(\mathcal A)$ of simple objects.

Using the 6-term exact sequences in cohomology, one may define a
biadditive form $\langle\cdot,\cdot\rangle$ on $K(\mathcal A)$ such that
for all objects $V$, $W$ of $\mathcal A$,
\begin{equation*}
\langle d(V),d(W)\rangle=\dim\Hom_{\mathcal A}(V,W)-\dim\Ext^1_{\mathcal A}(V,W).
\end{equation*}
This form is called the \textit{Euler form}.

\subsection{Hall numbers}\label{ss:HallNumbers}
In the remainder of Section~\ref{se:HallAlg}, the field $k$ will be the finite field $\Fq$ with $q$
elements.

Given an isomorphism class $\alpha\in\Iso(\mathcal A)$, we denote
the order of the automorphism group $\Aut_{\mathcal A}(M(\alpha))$ by $g_\alpha$.

Given three isomorphism classes $\alpha,\beta,\gamma\in\Iso(\mathcal A)$, we denote by
$\phi_{\alpha\gamma}^\beta$ the number of subobjects $X\subseteq M(\beta)$ such that
$X\in\gamma$ and $M(\beta)/X\in\alpha$. To be more precise, let $S(\alpha,\beta,\gamma)$ be
the set of pairs
\begin{equation*}
(f,g)\in\Hom_{\mathcal A}(M(\gamma),M(\beta))\times\Hom_{\mathcal A}(M(\beta),M(\alpha))
\end{equation*}
such that the sequence
\begin{equation*}
0\longrightarrow M(\gamma)\,{\buildrel f\over\longrightarrow}\,M(\beta)\,{\buildrel
g\over\longrightarrow}\,M(\gamma)\longrightarrow0
\end{equation*}
is exact. The group $\Aut_{\mathcal A}(M(\alpha))\times\Aut_{\mathcal A}(M(\gamma))$
acts freely on $S(\alpha,\beta,\gamma)$, and $\phi_{\alpha\gamma}^\beta$ is by
definition the cardinality of the quotient space
$S(\alpha,\beta,\gamma)/(\Aut_{\mathcal A}(M(\alpha))\times\Aut_{\mathcal A}(M(\gamma)))$.

The integer $\phi_{\alpha\gamma}^\beta$ is called a \textit{Hall number}. It is
zero if $d(\beta)\neq d(\alpha)+d(\gamma)$. Hall numbers have
the following properties.
\begin{proposition}\label{pr:PropHallNum}
If $\alpha,\beta,\gamma,\delta\in\Iso(\mathcal A)$ are isomorphism classes, then
\begin{enumerate}
\item\label{it:PrPHNa}there are only finitely many isomorphism classes $\lambda$ such that
$\phi_{\alpha\gamma}^\lambda\neq0$;
\item\label{it:PrPHNb}if (H5) holds, there are only finitely many pairs $(\rho,\sigma)\in
\Iso(\mathcal A)^2$ such that $\phi_{\rho\sigma}^\beta\neq0$;
\item\label{it:PrPHNc}$\phi_{\alpha0}^\beta=\delta_{\alpha\beta}$ and $\phi_{0\gamma}^\beta=
\delta_{\beta\gamma}$ (Kronecker symbols);
\item\label{it:PrPHNd}$\displaystyle\sum_{\lambda\in\Iso(\mathcal A)}\phi_{\alpha\beta}^\lambda
\phi_{\lambda\gamma}^\delta=\sum_{\lambda\in\Iso(\mathcal A)}\phi_{\alpha\lambda}^\delta
\phi_{\beta\gamma}^\lambda$;
\item\label{it:PrPHNe}$q^{\dim\Hom_{\mathcal A}(M(\alpha),M(\gamma))}\;\phi_{\alpha
\gamma}^\beta\;g_\alpha g_\gamma/g_\beta$ is an integer;
\item\label{it:PrPHNf}$\displaystyle\sum_{\lambda\in\Iso(\mathcal A)}\phi_{\alpha\gamma}^\lambda\;
g_\alpha g_\gamma/g_\lambda=q^{-\langle d(\alpha),d(\gamma)\rangle}$;
\item\label{it:PrPHNg}if $M(\alpha)$ and $M(\gamma)$ are indecomposable objects and $M(\beta)$ is
a decomposable object, then $q-1$ divides $\phi_{\alpha\gamma}^\beta-\phi_{\gamma\alpha}^\beta$;
\item\label{it:PrPHNh}the following formula holds:
\begin{equation*}
g_\alpha g_\beta g_\gamma g_\delta\sum_{\lambda\in\Iso(\mathcal A)}\phi_{\alpha\beta}^\lambda
\,\phi_{\gamma\delta}^\lambda/g_\lambda=\sum_{\rho,\rho',\sigma,\sigma'\in\Iso(\mathcal A)}
q^{-\langle d(\rho),d(\sigma')\rangle}\,\phi_{\rho\rho'}^\alpha\,
\phi_{\sigma\sigma'}^\beta\,\phi_{\rho\sigma}^\gamma\,\phi_{\rho'\sigma'}^\delta\,g_\rho
g_{\rho'} g_\sigma g_{\sigma'}.
\end{equation*}
\end{enumerate}
\end{proposition}
In the above statement, the sums in Items \ref{it:PrPHNd}, \ref{it:PrPHNf} and \ref{it:PrPHNh}
involve a finite number of non-zero terms.
\begin{proof}
Assertion~\ref{it:PrPHNa} holds because the extension group $\Ext^1_{\mathcal A}(M(\alpha),M(\gamma))$
is a finite set. Conditions~(H3) and (H5) imply that the map $d:\Iso(\mathcal A)\to K(\mathcal A)$
has finite fibers; Assertion~\ref{it:PrPHNb} follows from this fact. Assertion~\ref{it:PrPHNc}
is trivial. Assertion~\ref{it:PrPHNd} is Proposition~1 in~\cite{Ringel2}.
Items~\ref{it:PrPHNe} and \ref{it:PrPHNf} are consequences of Proposition~I.3.4 in~\cite{Ringel5}
(see also~\cite{Macdonald}, p.~221). To prove~\ref{it:PrPHNg}, it suffices to follow the proof of
Proposition~1 in~\cite{Ringel3}. Finally, to prove~\ref{it:PrPHNh}, one can adapt the proof of
Theorem~2 in~\cite{Green1} to the present framework.
\end{proof}

\subsection{The Hall algebra and the Ringel-Green bialgebra}\label{ss:HallAlgebra}
We will use $\widetilde{\mathbf Z}=\mathbf Z[v,v^{-1}]/(v^2-q)$ as the ground ring.
Let $H(\mathcal A)$ be the free $\widetilde{\mathbf Z}$-module on the symbols $\alpha$, where $\alpha$
runs over $\Iso(\mathcal A)$. The multiplication
\begin{equation*}
\alpha\cdot\gamma=\sum_{\beta\in\Iso(\mathcal A)}\phi_{\alpha\gamma}^\beta\;\beta
\end{equation*}
endows $H(\mathcal A)$ with the structure of an associative $\widetilde{\mathbf Z}$-algebra with
unit given by $0$. This follows from Items \ref{it:PrPHNa}, \ref{it:PrPHNc}, and \ref{it:PrPHNd} of
Proposition~\ref{pr:PropHallNum}. The algebra $H(\mathcal A)$ is called the \textit{Hall algebra}
of the category $\mathcal A$. It is graded by the group $K(\mathcal A)$, the
symbol $\alpha$ being homogeneous of degree~$d(\alpha)$.

Ringel~\cite{Ringel1} observed that the Euler form of Section~\ref{ss:GrotEuler} allows to equip
$H(\mathcal A)$ with another multiplication~$*$, which is defined on the basis vectors of
$H(\mathcal A)$ by
\begin{equation*}
\alpha*\gamma=v^{\langle d(\alpha),d(\gamma)\rangle}\;\alpha\cdot\gamma.
\end{equation*}
Since the Euler form is biadditive, this law endows $H(\mathcal A)$ with another structure
of an associative $\widetilde{\mathbf Z}$-algebra, called the \textit{Ringel algebra} of $\mathcal A$.

If Condition (H5) holds, then one can define a coproduct $\Delta:
H(\mathcal A)\to H(\mathcal A)\otimes_{\widetilde{\mathbf Z}}H(\mathcal A)$
and a co\"unit $\varepsilon:H(\mathcal A)\to\widetilde{\mathbf Z}$ by
\begin{equation*}
\Delta(\beta)=\sum_{\alpha,\gamma\in\Iso(\mathcal A)}v^{\langle d(\alpha),
d(\gamma)\rangle}\;{g_\alpha g_\gamma\over g_\beta}\;\phi_{\alpha\gamma}^\beta\
(\alpha\otimes\gamma)\qquad\text{and}\qquad\varepsilon(\beta)=\delta_{\beta0},
\end{equation*}
for all $\beta\in\Iso(\mathcal A)$. In this way, $H(\mathcal A)$ becomes a
$\widetilde{\mathbf Z}$-coalgebra in view of Properties
\ref{it:PrPHNb}--\ref{it:PrPHNe} of Proposition~\ref{pr:PropHallNum}.
Property~\ref{it:PrPHNh} of Proposition~\ref{pr:PropHallNum}
implies that $\Delta$ is an homomorphism of algebras when one equips $H(\mathcal A)
\otimes_{\widetilde{\mathbf Z}}H(\mathcal A)$ with the following twisted product:
\begin{equation*}
(\alpha\otimes\beta)*(\gamma\otimes\delta)=v^{\langle d(\beta),d(\gamma)
\rangle+\langle d(\gamma),d(\beta)\rangle}\;(\alpha*\gamma)\otimes(\beta*\delta),
\end{equation*}
where $\alpha,\beta,\gamma,\delta\in\Iso(\mathcal A)$. Endowed with the Ringel product
$*$, the coproduct $\Delta$, and the co\"unit~$\varepsilon$, the $\widetilde{\mathbf Z}$-module
$H(\mathcal A)$ is called the \textit{twisted Ringel-Green bialgebra}.

\begin{other}{Remark}
Suppose $\mathcal A$ satisfies (H5). In order to turn the twisted Ringel-Green
bialgebra into an actual bialgebra, Jie Xiao \cite{Xiao} and Kapranov \cite{Kapranov}
proceed as follows. They observe that the group $K(\mathcal A)$ acts by automorphisms
on the Ringel algebra $H(\mathcal A)$ by
\begin{equation*}
x\cdot\alpha=v^{-\langle x,d(\alpha)\rangle-\langle d(\alpha),x\rangle}\;\alpha,
\end{equation*}
where $x\in K(\mathcal A)$ and $\alpha\in\Iso(\mathcal A)$. Using this action, they
form the twisted group algebra $B(\mathcal A)$ of $K(\mathcal A)$ over $H(\mathcal A)$.
The algebra $B(\mathcal A)$ is a free $\widetilde{\mathbf Z}$-module
with basis the set of symbols $K_x\otimes\alpha$, where $x\in K(\mathcal A)$ and
$\alpha\in\Iso(\mathcal A)$. The multiplication on $B(\mathcal A)$ is given by
\begin{equation*}
(K_x\otimes\alpha)(K_y\otimes\beta)=v^{\langle y,d(\alpha)\rangle+\langle
d(\alpha),y\rangle}\ K_{x+y}\otimes(\alpha*\beta),
\end{equation*}
for all $x,y\in K(\mathcal A)$ and $\alpha,\beta\in\Iso(\mathcal A)$.

Endowed with the coproduct $\Delta:B(\mathcal A)\to B(\mathcal A)\otimes_{\widetilde{\mathbf Z}}
B(\mathcal A)$ and the co\"unit $\varepsilon:B(\mathcal A)\to\widetilde{\mathbf Z}$ given~by
\begin{align*}
\Delta(K_x\otimes\beta)&\,=\,\sum_{\alpha,\gamma\in\Iso(\mathcal A)}v^{\langle d(\alpha),
d(\gamma)\rangle}\;{g_\alpha g_\gamma\over g_\beta}\;\phi_{\alpha\gamma}^\beta\ 
(K_x\otimes\alpha)\otimes(K_{x+d(\alpha)}\otimes\gamma)\\
\intertext{and}
\varepsilon(K_x\otimes\beta)&\,=\;\delta_{\beta0},
\end{align*}
$B(\mathcal A)$ becomes a $\widetilde{\mathbf Z}$-bialgebra.
Jie Xiao and Kapranov further show that $B(\mathcal A)$ has an antipode.
\end{other}

\section{Coherent sheaves over the projective line}\label{se:CohSheavesP1}
Let $k$ be a field. In this section, we investigate the category $\mathcal A$ of coherent
sheaves over the projective line $\Pone(k)$. We determine the indecomposable objects of
$\mathcal A$ and study certain extensions between them. This information will be used in
Section~\ref{se:HallAlgCohShP1} to determine structure constants of the Hall algebra
$H(\mathcal A)$ when $k$ is a finite field.

\subsection{Generalities on coherent sheaves on $\Pone(k)$}\label{ss:CohShP1}
We put homogeneous coordinates $(t:u)$ on $\Pone(k)$. The two affine open subsets
\begin{equation*}
U'=\{(t:u)\mid t\neq0\}\quad\text{and}\quad U''=\{(t:u)\mid u\neq0\}
\end{equation*}
cover $\Pone(k)$, and the formulae $z=u/t$ and $z^{-1}=t/u$ define coordinates on
$U'$ and $U''$ respectively. The rings $k[z]$ and $k[z^{-1}]$ are the respective
rings of regular functions on $U'$ and~$U''$.

A closed point $x$ of $\Pone(k)$ is the zero locus of an irreducible homogeneous polynomial
$P\in k[T,U]$. If $P$ is proportional to the polynomial $T$, then the closed point is the point
at infinity $\infty$. If $P$ is not proportional to $T$, then $x$ can be viewed as
the zero locus in $U'$ of the irreducible polynomial $P(1,z)\in k[z]$. In any case, $x$ determines
$P$ up to a non-zero scalar, and the degree $\deg x$ of $x$ is defined as the degree of $P$.

We will use the following convention: if $A$ is a commutative ring, $M$
an $A$-module and $z$ an element of $A$, then $M_z$ denotes the localized
$A$-module obtained from $M$ by inverting~$z$. An analogous notation will be used for morphisms.

We are now ready to define the category $\mathcal A$ of \textit{coherent sheaves}
on $\Pone(k)$. An object of $\mathcal A$ is a triple $(M',M'',\varphi)$,
where $M'$ is a finitely generated $k[z]$-module, $M''$ is a finitely generated
$k[z^{-1}]$-module, and $\varphi:M'_z\to M''_{z^{-1}}$ is an isomorphism of
$k[z,z^{-1}]$-modules. A morphism in $\mathcal A$ from the coherent sheaf
$(M',M'',\varphi)$ to the coherent sheaf $(N',N'',\psi)$ is a pair of maps $(f',f'')$,
where $f':M'\to N'$ is a $k[z]$-linear map and $f'':M''\to N''$ is a $k[z^{-1}]$-linear
map such that $\psi\circ f'_z=f''_{z^{-1}}\circ\varphi$.
One also defines in an obvious way the notions of direct sums and exact sequences
in $\mathcal A$, so that $\mathcal A$ becomes an abelian $k$-linear category.
This definition of $\mathcal A$ is equivalent to the standard geometric
definition that can be found, for instance, in Section~II.5 of~\cite{Hartshorne}.

A coherent sheaf $(M',M'',\varphi)$ is said to be \textit{locally free} if $M'$ and
$M''$ are free modules over $k[z]$ and $k[z^{-1}]$ respectively. The full subcategory
of $\mathcal A$ consisting of all locally free sheaves will be denoted by $\mathcal A_\lf$.

For any $n\in\mathbf Z$, we construct a locally free coherent sheaf $(M',M'',\varphi)$ by
letting $M'=k[z]$, $M''=k[z^{-1}]$, and $\varphi:k[z,z^{-1}]\to k[z,z^{-1}]$ be the
multiplication by $z^{-n}$. As usual, this sheaf will be denoted by $\mathcal O(n)$.
(The structure sheaf $\mathcal O_{\Pone(k)}$ of the projective line is the sheaf~$\mathcal O(0)$.)
For any $m,n\in\mathbf Z$, the space of homogeneous polynomials $F\in k[T,U]$ of degree $n-m$
is naturally isomorphic to the homomorphism space $\Hom_{\mathcal A}(\mathcal O(m),\mathcal O(n))$:
one associates to $F$ the pair of maps $(f',f'')$, where $f':k[z]\to k[z]$
is the multiplication by $F(1,z)$ and $f'':k[z^{-1}]\to k[z^{-1}]$ is the
multiplication by~$F(z^{-1},1)$.

A coherent sheaf $(M',M'',\varphi)$ is called a \textit{torsion sheaf} if $M'$
is a torsion $k[z]$-module, which is equivalent to a similar requirement for $M'_z$,
$M''_{z^{-1}}$, or $M''$. The full subcategory of $\mathcal A$ consisting of all torsion
sheaves will be denoted by $\mathcal A_\tor$.

Given an irreducible homogeneous polynomial $P\in k[T,U]$ of degree~$d$ and an integer
$r\geq 1$, the $r$-th power polynomial $P^r$ defines a morphism from $\mathcal O(-rd)$ to
$\mathcal O(0)$. The cokernel is the torsion sheaf $(M',M'',\varphi)$, where
$M'=k[z]/(P(1,z)^r)$, $M''=k[z^{-1}]/(P(z^{-1},1)^r)$ and $\varphi$ is induced by the
identity of $k[z,z^{-1}]$. If $x$ is the closed point corresponding to
$P$, we denote this torsion sheaf by $\mathcal O_{r[x]}$.

\begin{proposition}\label{pr:IndecCatA}
\begin{enumerate}
\item\label{it:PrICAa}The category $\mathcal A$ is $k$-linear, abelian, and satisfies
Conditions~(H1)--(H4) of Section~\ref{ss:InitData}. The subcategories $\mathcal A_\lf$
and $\mathcal A_\tor$ of $\mathcal A$ are closed under extensions. The category
$\mathcal A_\tor$ is $k$-linear, abelian, and satisfies Conditions~(H1)--(H5)
of Section~\ref{ss:InitData}.
\item\label{it:PrICAb}Every coherent sheaf $\mathcal F$ can be written as a direct sum
$\mathcal F_0\oplus\mathcal F_1$, where $\mathcal F_0$ is a torsion sheaf and $\mathcal F_1$
is locally free. The isomorphism classes of $\mathcal F_0$ and $\mathcal F_1$ are
determined by the isomorphism class of $\mathcal F$.
\item\label{it:PrICAc}Up to isomorphism, the indecomposable objects in $\mathcal A$ are
the locally free sheaves $\mathcal O(n)$, where $n\in\mathbf Z$, and the torsion sheaves
$\mathcal O_{r[x]}$, where $r$ is an integer $\geq1$ and $x$ is a closed point of $\Pone(k)$.
\end{enumerate}
\end{proposition}
\begin{proof}
\begin{enumerate}
\item It is clear that $\mathcal A$ is abelian and satisfies Condition~(H1). Propositions~III.2.2,
III.6.3~(c), III.6.4, and III.6.7, Theorems~III.2.7 and III.5.2~(a), and Exercise~III.6.5
of~\cite{Hartshorne} imply that $\mathcal A$ satisfies Conditions~(H2)--(H4). It follows from
the definitions that $\mathcal A_\lf$ and $\mathcal A_\tor$ are closed under extensions
in $\mathcal A$, so that $\mathcal A_\lf$ and $\mathcal A_\tor$ also satisfy Conditions
(H1)--(H4). Since subobjects and quotients of torsion sheaves are torsion sheaves,
$\mathcal A_\tor$ is abelian. Finally, the simple fact that finitely generated torsion modules
over principal ideal domains have Jordan-H\"older series implies that $\mathcal A_\tor$
satisfies Condition~(H5).
\item Let $\mathcal F$ be a coherent sheaf. One can define in an obvious way the torsion
subsheaf $\tor(\mathcal F)$ of $\mathcal F$, and the quotient sheaf $\mathcal F/\tor(\mathcal F)$
is locally free. By Serre's vanishing theorem (Theorem~III.5.2~(b) in~\cite{Hartshorne}), the
extension group $\Ext^1_{\mathcal A}(\mathcal F/\tor(\mathcal F),\tor(\mathcal F))$ vanishes.
Thus the short exact sequence
\begin{equation*}
0\longrightarrow\tor(\mathcal F)\longrightarrow\mathcal F\longrightarrow\mathcal F/
\tor(\mathcal F)\longrightarrow0
\end{equation*}
splits and one gets the decomposition $\mathcal F\simeq\tor(\mathcal F)\oplus\mathcal F/\tor(\mathcal F)$.
Conversely, given a decomposition $\mathcal F=\mathcal F_0\oplus\mathcal F_1$ as in the
statement of Assertion~\ref{it:PrICAb}, one has $\mathcal F_0=\tor(\mathcal F)$ and $\mathcal F_1
\simeq\mathcal F/\mathcal F_0=\mathcal F/\tor(\mathcal F)$.
\item By Assertion~\ref{it:PrICAb}, an indecomposable coherent sheaf is either a torsion sheaf
or a locally free sheaf. The classification of torsion modules over the principal ideal domains
$k[z]$ and $k[z^{-1}]$ leads to the fact that the indecomposables torsion coherent sheaves
are the sheaves $\mathcal O_{r[x]}$, where $r\geq1$ and $x$ is a closed point of $\Pone(k)$.
On the other hand, a theorem of Grothendieck \cite{Grothendieck} asserts that any locally free
coherent sheaf is isomorphic to a direct sum $\mathcal O(n_1)\oplus\cdots\oplus\mathcal O(n_r)$
for some sequence $n_1\leq\cdots\leq n_r$ of uniquely determined integers. Thus the
sheaves $\mathcal O(n)$ are the indecomposable locally free coherent sheaves.
\end{enumerate}
\end{proof}

\subsection{The Grothendieck group and the Euler form}\label{ss:GrotCohShP1}
We define the rank and the degree of an indecomposable sheaf by
\begin{equation*}
\rk\mathcal O(n)=1,\quad\deg\mathcal O(n)=n,\quad\rk\mathcal O_{r[x]}=0,\quad\text{and}
\quad\deg\mathcal O_{r[x]}=r\,\deg x.
\end{equation*}
Since every coherent sheaf can be written in an essentially unique way as a sum of
indecomposable sheaves, we may extend additively the notions of rank and degree
to arbitrary coherent sheaves. Note that the torsion sheaves are the sheaves whose rank is $0$.

It is well-known from algebraic geometry that the rank and degree maps factor through the
Grothendieck group $K(\mathcal A)$, defining a morphism of abelian groups
$\upsilon:K(\mathcal A)\to\mathbf Z^2$ by
\begin{equation*}
d(\mathcal F)\mapsto(\rk\mathcal F,\deg\mathcal F).
\end{equation*}
\begin{proposition}\label{pr:GrotGroup}
\begin{enumerate}
\item\label{it:PrGGa}The homomorphism $\upsilon:K(\mathcal A)\to\mathbf Z^2$ is an isomorphism
of abelian groups.
\item\label{it:PrGGb}The Euler form on $K(\mathcal A)$ is given for all coherent sheaves
$\mathcal F$ and $\mathcal G$ by
\begin{equation}
\langle d(\mathcal F),d(\mathcal G)\rangle=\rk\mathcal F\;\rk\mathcal G+\rk\mathcal F
\;\deg\mathcal G-\deg\mathcal F\;\rk\mathcal G.\label{eq:RiemannRoch}
\end{equation}
\end{enumerate}
\end{proposition}
\begin{proof}
\begin{enumerate}
\item It is well known (for details, see Proposition~\ref{pr:ExtCS1}) that there are short exact
sequences of the form
\begin{equation*}
0\longrightarrow\mathcal O(m)\longrightarrow\mathcal O(a)\oplus
\mathcal O(m+n-a)\longrightarrow\mathcal O(n)\longrightarrow0
\end{equation*}
whenever $m\leq a\leq m+n-a\leq n$. Setting $m=0$, we get
\begin{equation*}
d(\mathcal O(n))=d(\mathcal O(a))+d(\mathcal O(n-a))-d(\mathcal O(0))
\end{equation*}
when $0\leq a\leq n$, and an easy induction shows that
\begin{equation}
d(\mathcal O(n))=n\,d(\mathcal O(1))+(1-n)\,d(\mathcal O(0))\label{eq:GenGrotGroup}
\end{equation}
for all integers $n\geq1$. A similar argument shows the validity of (\ref{eq:GenGrotGroup})
for $n\leq0$.

Now take a closed point $x$ and an integer $r\geq1$. By definition, the sheaf $\mathcal O_{r[x]}$
is the cokernel of an homomorphism in $\Hom_{\mathcal A}(\mathcal O(-r\deg x),
\mathcal O(0))$, so
\begin{equation*}
d(\mathcal O_{r[x]})=d(\mathcal O(0))-d(\mathcal O(-r\deg x))=r\deg x\;
\left(d(\mathcal O(1))-d(\mathcal O(0))\right).
\end{equation*}
The above discussion shows that $d(\mathcal O(0))$ and $d(\mathcal O(1))$
generate the group $K(\mathcal A)$. Since $\upsilon(d(\mathcal O(0)))=(1,0)$
and $\upsilon(d(\mathcal O(1)))=(1,1)$ form a basis of $\mathbf Z^2$,
Statement~\ref{it:PrGGa} holds.
\item Using the Riemann-Roch formula (Theorem~IV.1.3 in~\cite{Hartshorne}) and
standard results of sheaf cohomology (Propositions~II.5.12, III.6.3~(c), and III.6.7
in~\cite{Hartshorne}), we obtain
\begin{eqnarray*}
\langle d(\mathcal O(m)),d(\mathcal O(n))\rangle&=&
\dim\Hom_{\mathcal A}(\mathcal O(m),\mathcal O(n))
- \dim\Ext^1_{\mathcal A}(\mathcal O(m),\mathcal O(n))\\
&=&\dim\Hom_{\mathcal A}(\mathcal O_{\Pone(k)},\mathcal O(n-m))-
\dim\Ext^1_{\mathcal A}(\mathcal O_{\Pone(k)},\mathcal O(n-m))\\
&=&\dim H^0(\Pone(k),\mathcal O(n-m))-\dim H^1(\Pone(k),\mathcal O(n-m))\\
&=&1+\deg(\mathcal O(n-m))\\
&=&1+n-m\\
&=&\rk\mathcal O(m)\;\rk\mathcal O(n)\;+\;\rk\mathcal O(m)\;\deg\mathcal O(n)\;
-\;\deg\mathcal O(m)\;\rk\mathcal O(n)
\end{eqnarray*}
for all $m,n\in\mathbf Z$.
This proves~(\ref{eq:RiemannRoch}) when $\mathcal F$ and $\mathcal G$ are locally free sheaves
of rank~1. Since the classes of $\mathcal O(0)$ and $\mathcal O(1)$ generate $K(\mathcal A)$,
the general case follows by the biadditivity of both sides of~(\ref{eq:RiemannRoch}).
\end{enumerate}
\end{proof}

\subsection{Extensions of locally free sheaves}\label{ss:HomCohShP1}
Our first original result, presented below, describes the extensions between the indecomposable
locally free sheaves. Beforehand, let us record the following easy lemma, which is a consequence
of the description of the homomorphism spaces $\Hom_{\mathcal A}(\mathcal O(m),\mathcal O(n))$
recalled in Section~\ref{ss:CohShP1}.

\begin{lemma}\label{le:HomLFCohShP1}
For all $m,n\in\mathbf Z$,
\begin{enumerate}
\item\label{it:LeHLFCSPa}any non-zero element in $\Hom_{\mathcal A}(\mathcal O(m),\mathcal O(n))$ is
a monomorphism;
\item\label{it:LeHLFCSPb}as a $k$-algebra, $\End_{\mathcal A}(\mathcal O(n))\simeq k$;
\item\label{it:LeHLFCSPc}the $k$-vector space $\Hom_{\mathcal A}(\mathcal O(m),\mathcal O(n))$
has dimension $\max(0,n-m+1)$.
\end{enumerate}
\end{lemma}

We now analyse the short exact sequences of the form
\begin{equation*}
0\longrightarrow\mathcal O(m)\,{\buildrel f\over\longrightarrow}
\,\mathcal F\,{\buildrel g\over\longrightarrow}\,\mathcal O(n)
\longrightarrow0.
\end{equation*}
By Proposition~\ref{pr:IndecCatA}~\ref{it:PrICAa}, the coherent sheaf $\mathcal F$ is necessarily
locally free. Using rank considerations and Proposition~\ref{pr:IndecCatA}~\ref{it:PrICAc}, we
may assume, without any loss, that $\mathcal F$ is the sheaf $\mathcal O(p)\oplus\mathcal O(q)$
for some integers $p$ and $q\in\mathbf Z$.

\begin{proposition}\label{pr:ExtCS1}
Let $m$, $n$, $p$, $q$ be integers, and consider a sequence of the form
\begin{equation*}
0\longrightarrow\mathcal O(m)\,{\buildrel f\over\longrightarrow}\,\mathcal O(p)\oplus
\mathcal O(q)\,{\buildrel g\over\longrightarrow}\,\mathcal O(n)\longrightarrow0.
\end{equation*}
Let
\begin{eqnarray*}
h\in\Hom_{\mathcal A}(\mathcal O(m),\mathcal O(p)),\quad&\qquad&\quad
j\in\Hom_{\mathcal A}(\mathcal O(p),\mathcal O(n)),\\
i\in\Hom_{\mathcal A}(\mathcal O(m),\mathcal O(q)),\quad&&\quad
\ell\in\Hom_{\mathcal A}(\mathcal O(q),\mathcal O(n)),
\end{eqnarray*}
be defined by $f=h\oplus i$ and $g=j\oplus\ell$, and call $H$, $I$, $J$, $L$ the homogeneous polynomials
in $k[T,U]$ representing $h$, $i$, $j$, $\ell$, respectively. Then the sequence is a non-split short
exact sequence if and only if the following three conditions are satisfied:
\begin{description}
\item[(a)]$m<\min(p,q)$, \ $\max(p,q)<n$, \ $p+q=m+n$.
\item[(b)]$J$ and $L$ are coprime polynomials.
\item[(c)]There is a non-zero scalar $E$ such that $H=EL$ and $I=-EJ$.
\end{description}
\end{proposition}
\begin{proof}
We first prove that Conditions (a), (b), and (c) are necessary. Suppose that the sequence is
exact and non-split. If one of the homomorphism $h$ or $i$ were the zero arrow,
then the other one would be an isomorphism since the cokernel of $f$
is indecomposable, and the sequence would split. Similarly, neither $j$ nor $\ell$ can vanish.
The four maps $h$, $i$, $j$, and $\ell$ are thus non-zero, and none of them is an isomorphism.
In view of Lemma~\ref{le:HomLFCohShP1}~\ref{it:LeHLFCSPb} and \ref{it:LeHLFCSPc}, it follows
that $m<p,q<n$. For degree reasons we also have $m+n=p+q$. Therefore Condition~(a) holds.

Let us turn to Condition~(b). In the unique factorization domain $k[T,U]$, one may
consider a g.c.d.\ $D$ of the polynomials $J$ and $L$. Since every irreducible factor of a
homogeneous polynomial is itself homogeneous (by uniqueness of the factorization), the
polynomial $D$ is homogeneous and define a homomorphism $d\in\Hom_{\mathcal A}(\mathcal
O(n-\deg D),\mathcal O(n))$. We get a factorization
\begin{equation*}
\begin{matrix}
\mathcal O(p)\oplus\mathcal O(q)&\buildrel g\over\longlongrightarrow&\mathcal O(n)\\[2pt]
\qquad\searrow&&\nearrow{\scriptstyle d}\qquad\\&\mathcal O(n-\deg D)&\end{matrix}\ .
\end{equation*}
Since $g$ is surjective, so must be $d$. The non-zero morphism $d$ being injective
by Lemma~\ref{le:HomLFCohShP1}~\ref{it:LeHLFCSPa}, it is an isomorphism,
which implies that $\deg D=0$. Thus $J$ and $L$ are coprime, and Condition~(b) holds.

Finally, the equality $g\circ f=0$ implies that $HJ+IL=0$. Condition~(b) and Gauss's lemma
then imply the existence of a non-zero homogeneous polynomial $E\in k[T,U]$ such that
$H=EL$ and $I=-EJ$. Since
\begin{equation*}
2\deg E=\deg H+\deg I-\deg L-\deg J=(p-m)+(q-m)-(n-q)-(n-p)=0,
\end{equation*}
$E$ is a constant polynomial, which proves Condition~(c).

In order to prove the converse statement, we now assume that Conditions~(a), (b), and (c)
are fulfilled. Over the affine subset $U'$, our sequence of sheaves reads
\begin{equation*}
\begin{matrix}0\longrightarrow&\mathcal O(m)(U')&\buildrel f_{U'}=h_{U'}\oplus i_{U'}
\over\longlongrightarrow&\mathcal O(p)(U')\oplus\mathcal O(q)(U')&\buildrel g_{U'}=
j_{U'}\oplus\ell_{U'}\over\longlongrightarrow&\mathcal O(n)(U')&\longrightarrow0 ,\\
&\|&&\|&&\|&\cr&k[z]&&k[z]\oplus k[z]&&k[z]&\end{matrix}
\end{equation*}
where the maps $h_{U'}$, $i_{U'},\ldots$ are the multiplications by $H(1,z)$,
$I(1,z),\ldots$ respectively. Condition~(b) implies that the polynomials $J(1,z)$ and
$L(1,z)$ are coprime, which ensures by Bezout's lemma that the $k[z]$-linear map $g_{U'}$
is surjective. An analogous simple reasoning based on Gauss's lemma shows that Conditions~(b)
and (c) imply that $\ker g_{U'}=\im f_{U'}$. Thus our sequence of sheaves is exact over
the open subset $U'$. A similar argument can be used over $U''$, and we conclude
that our sequence of sheaves is exact.
\end{proof}

\begin{corollary}\label{co:VanishExt}
If $m,n\in\mathbf Z$ are integers satisfying $n\leq m+1$, then the extension group
$\Ext^1_{\mathcal A}(\mathcal O(n),\mathcal O(m))$ vanishes.
\end{corollary}

\subsection{Torsion sheaves}\label{ss:TorsionSheaves}
Torsion sheaves will play an important r\^ole in Section~\ref{se:HallAlgCohShP1}. Although
they are simpler objects than locally free sheaves, their precise description is rather
technical. In this section, we set some further notation and we establish some basic facts.

According to standard terminology, a partition is a non-increasing
sequence of non-negative integers with only finitely many non-zero terms:
$\lambda=(\lambda_1\geq\lambda_2\geq\cdots)$ with $\lambda_i=0$ for $i$ big enough.
The length of $\lambda$ is the smallest integer $\ell\geq0$ such that $\lambda_{\ell+1}=0$,
and the weight $|\lambda|$ of $\lambda$ is the sum of the non-zero integers $\lambda_i$.
We also put \vspace{2pt}$n(\lambda)=\sum_{i\geq1}(i-1)\lambda_i$, as in~\cite{Macdonald}.
The empty partition is the partition with no non-zero part; the partition with $r$ non-zero parts,
all equal to $1$, is denoted by $(1^r)$; the partition with one non-zero part, equal to $r$, is
denoted by $(r)$.

For any closed point $x$ and any partition $\lambda=(\lambda_1\geq\lambda_2\geq\cdots)$ of
length $\ell$, we define the torsion sheaf
\begin{equation*}
\mathcal O_{\lambda[x]}=\mathcal O_{\lambda_1[x]}\oplus\cdots\oplus\mathcal O_{\lambda_\ell[x]}.
\end{equation*}
For instance, $\mathcal O_{(1^r)[x]}=\left(\mathcal O_{[x]}\right)^{\oplus r}$ and $\mathcal O_{(r)[x]}=
\mathcal O_{r[x]}$. By Proposition~\ref{pr:IndecCatA}~\ref{it:PrICAc}, for any torsion sheaf $\mathcal F$,
there is a finite family $(x_i)_{1\leq i\leq t}$ of distinct closed points of $\Pone(k)$ and a finite family
$(\lambda^{(i)})_{1\leq i\leq t}$ of non-empty partitions such that
\begin{equation*}
\mathcal F\simeq\mathcal O_{\lambda^{(1)}[x_1]}\oplus\cdots\oplus\mathcal O_{\lambda^{(t)}[x_t]}.
\end{equation*}
The set $\{x_1,\ldots,x_t\}$, uniquely determined by the sheaf $\mathcal F$, is called the
\textit{support} of $\mathcal F$. We will denote by $\mathcal A_{\{x\}}$ the full subcategory
of $\mathcal A$ consisting of all torsion sheaves with support included in $\{x\}$.

\begin{lemma}\label{le:HomTCohShP1}
\begin{enumerate}
\item\label{it:LeHTCSPa}If $\mathcal F$ and $\mathcal G$ are torsion sheaves with disjoint supports,
then $\Hom_{\mathcal A}(\mathcal F,\mathcal G)=\Ext^1_{\mathcal A}(\mathcal F,\mathcal G)=0$.
\item\label{it:LeHTCSPb}Let $x$ be a closed point of $\Pone(k)$. The category $\mathcal A_{\{x\}}$
is $k$-linear, abelian, and satisfies Conditions~(H1)--(H5) of Section~\ref{ss:InitData}.
\item\label{it:LeHTCSPc}The category $\mathcal A_\tor$ is the direct sum of the subcategories
$\mathcal A_{\{x\}}$, where $x$ runs over the set of all closed points of $\Pone(k)$.
\item\label{it:LeHTCSPd}If $\mathcal F$ is a locally free sheaf and $\mathcal G$ a torsion sheaf, then
$\Hom_{\mathcal A}(\mathcal G,\mathcal F)=\Ext^1_{\mathcal A}(\mathcal F,\mathcal G)=0$.
\item\label{it:LeHTCSPe}For any closed point $x$, any partition $\lambda$, and any
$n\in\mathbf Z$, the $k$-vector space\break$\Hom_{\mathcal A}(\mathcal O(n),\mathcal O_{\lambda[x]})$ has
dimension $|\lambda|\deg x$.
\end{enumerate}
\end{lemma}
\begin{proof}
Assertion~\ref{it:LeHTCSPa} is a consequence of the structure theorem for finitely generated
torsion modules over a principal ideal domain. This assertion implies that
the class of sheaves with support in some fixed finite set $\{x_1,\ldots,x_t\}$ of closed points
of $\Pone(k)$ is closed under the operations of taking subobjects, quotients, and extensions.
Assertion~\ref{it:LeHTCSPb} follows therefore from Assertion~\ref{it:LeHTCSPa} and
Proposition~\ref{pr:IndecCatA}~\ref{it:PrICAa}. Since every torsion sheaf can be uniquely
written as a direct sum of subsheaves belonging to subcategories $\mathcal A_{\{x\}}$,
Assertion~\ref{it:LeHTCSPc} follows from Assertion~\ref{it:LeHTCSPa}.

Now let $\mathcal F$ be a locally free sheaf and $\mathcal G$ be a torsion sheaf. The
vanishing of $\Hom_{\mathcal A}(\mathcal G,\mathcal F)$ is a direct consequence of the
definitions in Section~\ref{ss:CohShP1}, while that of $\Ext^1_{\mathcal A}(\mathcal F,
\mathcal G)$ follows from Serre's vanishing theorem. Thus Assertion~\ref{it:LeHTCSPd} holds.

Finally, let $n\in\mathbf Z$, $\lambda$ be a partition, and $x$ be a closed point of
$\Pone(k)$. The sheaf $\mathcal O_{\lambda[x]}$ has rank $0$ and degree $|\lambda|\deg x$,
and the extension group $\Ext^1_{\mathcal A}(\mathcal O(n),\mathcal O_{\lambda[x]})$
vanishes. Using Proposition~\ref{pr:GrotGroup}~\ref{it:PrGGb}, we therefore get
\begin{align*}
\dim\Hom_{\mathcal A}(\mathcal O(n),\mathcal O_{\lambda[x]})
&=\langle\mathcal O(n),\mathcal O_{\lambda[x]}\rangle\\
&=\rk\mathcal O(n)\;\rk\mathcal O_{\lambda[x]}\;+\;\rk\mathcal O(n)\;\deg\mathcal O_{\lambda[x]}\;-\;
\deg\mathcal O(n)\;\rk\mathcal O_{\lambda[x]}\\
&=\rk\mathcal O(n)\;\deg\mathcal O_{\lambda[x]}\\
&=|\lambda|\deg x.
\end{align*}
\end{proof}

Further properties of the categories $\mathcal A_{\{x\}}$ (and of their Hall algebras when
$k$ is a finite field) required in Section~\ref{se:HallAlgCohShP1} rely on an isomorphism
of categories, which we now describe. Let $x$ be a fixed closed point of $\Pone(k)$, defined
by an irreducible homogeneous polynomial $P\in k[T,U]$. If $x$ belongs to the affine open set $U'$
(respectively to $U''$), one defines the local ring $\mathcal O_{\Pone(k),x}$ of rational functions
regular near $x$ as the localization of $k[z]$ at the prime ideal generated by $P(1,z)$
(respectively as the localization of $k[z^{-1}]$ at the prime ideal generated by $P(z^{-1},1)$).
If $x$ belongs to $U'$ and to $U''$, then both definitions of $\mathcal O_{\Pone(k),x}$ yield
isomorphic rings. We denote by $k_x$ the residue field of $\mathcal O_{\Pone(k),x}$; the field
$k_x$ is an extension of $k$ of degree $\deg x$. We choose a generator $\pi_x$ of the maximal
ideal of $\mathcal O_{\Pone(k),x}$ (``uniformizer''). We finally denote by 
$\catmodf{\mathcal O_{\Pone(k),x}}$ the category of $\mathcal O_{\Pone(k),x}$-modules of
finite length, that is, of modules that are finitely generated and annihilated by some power
of the uniformizer $\pi_x$. The following statement is clear from the definitions.

\begin{proposition}\label{pr:IsomCatA_x}
Let $x$ be a closed point of $\Pone(k)$. There is an exact additive $k$-linear functor from
$\catmodf{\mathcal O_{\Pone(k),x}}$ to $\mathcal A_{\{x\}}$, which yields an isomorphism
of categories, and which sends the $\mathcal O_{\Pone(k),x}$-module $\mathcal O_{\Pone(k),x}/
(\pi_x{}^r)$ to the sheaf $\mathcal O_{r[x]}$.
\end{proposition}

\subsection{Mixed extensions}\label{ss:MixedExt}
We now investigate the extensions of certain torsion sheaves by indecomposable locally free sheaves.
We begin with the following lemma.
\begin{lemma}\label{le:HomMiCohShP1}
Let $x$ be a closed point of $\Pone(k)$, corresponding to an irreducible homogeneous polynomial
$P\in k[T,U]$. Let $m,n\in\mathbf Z$, let $F\in k[T,U]$ be a non-zero homogeneous polynomial
of degree $n-m$, and let $f\in\Hom_{\mathcal A}(\mathcal O(m),\mathcal O(n))$ be the morphism
represented by $F$. If the support of the cokernel of $f$ is included in $\{x\}$, then there
exists an integer $r\geq1$ such that $F=P^r$, up to a non-zero scalar, and one has $\coker f\simeq
\mathcal O_{r[x]}$.
\end{lemma}
\begin{proof}
The polynomial $F$ represents a morphism $\widetilde f:\mathcal O(m-n)\to\mathcal O(0)$.
By unique factorization, and up to a non-zero scalar, we can write $F=P_1^{r_1}\cdots P_t^{r_t}$,
for some positive integers $r_1,\ldots,r_t$ and some distinct irreducible homogeneous polynomials
$P_1\ldots,P_t$. For each $1\leq i\leq t$, let $x_i$ be the closed point of $\Pone(k)$ corresponding
to $P_i$. The homogeneous polynomial $P_i^{r_i}$ defines an element of $\Hom_{\mathcal A}(\mathcal O
(-r_i\deg P_i),\mathcal O(0))$ whose cokernel is $\mathcal O_{r_i[x_i]}$, whence a canonical morphism
$g_i:\mathcal O(0)\to\mathcal O_{r_i[x_i]}$. A direct application of the Chinese remainder
theorem implies that the sequence
\begin{equation*}
0\longrightarrow\mathcal O(m-n)\,{\buildrel\widetilde f\over\longrightarrow}\,\mathcal O(0)\,{\buildrel
\oplus g_i\over\longrightarrow}\,\oplus_{i=1}^t\mathcal O_{r_i[x_i]}\longrightarrow0
\end{equation*}
is exact over the affine subsets $U'$ and $U''$, hence is exact. Taking the tensor product
with the locally free sheaf $\mathcal O(n)$, we get an exact sequence
\begin{equation*}
0\longrightarrow\mathcal O(m)\,{\buildrel f\over\longrightarrow}\,\mathcal O(n)
\longrightarrow\oplus_{i=1}^t\mathcal O_{r_i[x_i]}\longrightarrow0.
\end{equation*}
Therefore, the cokernel of $f$ is isomorphic to
$\mathcal O_{r_1[x_1]}\oplus\cdots\oplus\mathcal O_{r_t[x_t]}$. If the support of $\coker f$ is
included in $\{x\}$, then $t=1$ and $P_1=P$, up to a non-zero scalar, which entails the lemma.
\end{proof}

\begin{proposition}\label{pr:ExtCS2}
Given a closed point $x$ and an integer $r\geq1$, let
\begin{equation*}
0\longrightarrow\mathcal O(m)\,{\buildrel f\over\longrightarrow}\,\mathcal F\,{\buildrel
g\over\longrightarrow}\,\mathcal O_{(1^r)[x]}\longrightarrow0
\end{equation*}
be a non-split short exact sequence of coherent sheaves. Then the middle term $\mathcal F$ is
isomorphic to $\mathcal O_{(1^{r-1})[x]}\oplus\mathcal O(m+\deg x)$. If we write $f=h\oplus i$
in this decomposition, then $h=0$ and $\coker i\simeq\mathcal O_{[x]}$.
\end{proposition}
\begin{proof}
We can write $\mathcal F$ as the direct sum of its torsion subsheaf $\mathcal F_0=\tor(\mathcal F)$
and a locally free subsheaf $\mathcal F_1\simeq\mathcal F/\tor(\mathcal F)$ of rank~$1$. Write the
maps $f$ and $g$ as $h\oplus i$ and $j\oplus\ell$ in the decomposition $\mathcal F=
\mathcal F_0\oplus\mathcal F_1$. The morphism $h$ cannot be injective, so $i$ cannot be zero,
so $i$ is injective (Lemma~\ref{le:HomLFCohShP1}~\ref{it:LeHLFCSPa}), and it follows that $j$ is
injective. Thus $\mathcal F_0$ must be isomorphic to a subobject of
$\mathcal O_{(1^r)[x]}$.

Under the isomorphism of categories described in Proposition~\ref{pr:IsomCatA_x}, the sheaf
$\mathcal O_{(1^r)[x]}$ corresponds to the elementary \vspace{2pt}$\mathcal O_{\Pone(k),x}$-module
$\bigl(\mathcal O_{\Pone(k),x}/(\pi_x)\bigr)^{\oplus r}$, hence to a vector space of
dimension $r$ over the residue field $k_x$. This shows that $\mathcal F_0$ is isomorphic to
$\mathcal O_{(1^s)[x]}$ for some $s\leq r$. In the same way, and using Lemma~\ref{le:HomMiCohShP1},
we see that the image of $\ell$ is either $0$ or isomorphic to~$\mathcal O_{[x]}$.

Now if the sequence is not split, then $j$ is not an isomorphism, which rules out the
case $s=r$. The surjectivity of $g$ then requires that $s=r-1$, that $\im\ell\simeq
\mathcal O_{[x]}$, and that $\mathcal O_{(1^r)[x]}=\im j\oplus\im\ell$. The equality
$g\circ f=0$ then splits into the two equalities $j\circ h=0$ and $\ell\circ i=0$.
Since $j$ is injective, we get $h=0$ and thus $\coker i\simeq\mathcal O_{[x]}$.
Finally we compute
\begin{equation*}
\deg\mathcal F_1=\deg\mathcal O(m)+\deg\mathcal O_{(1^r)[x]}-\deg\mathcal F_0=
m+r\deg x-s\deg x=m+\deg x,
\end{equation*}
which shows that $\mathcal F_1\simeq\mathcal O(m+\deg x)$.
\end{proof}

\section{The Hall algebra of $\Coh(\Pone(\Fq))$}\label{se:HallAlgCohShP1}
From now on, $k$ is the finite field $\Fq$ with $q$ elements, while $\mathcal A$ still
stands for the category of coherent sheaves over $\Pone(k)$. In this section, we
describe the Ringel algebra $H(\mathcal A)$. Relying on results of Section~\ref{se:CohSheavesP1},
we first compute some Hall numbers. Using results explained in Chapters~II and III of~\cite{Macdonald},
we next investigate more closely the Ringel algebras $H(\mathcal A_{\{x\}})$ and
$H(\mathcal A_\tor)$, which will be viewed as subalgebras of $H(\mathcal A)$. This allows us
eventually to define a certain subalgebra $B$ of $H(\mathcal A)$, which will turn out
in Section~\ref{se:Uqsl2hat} to be related to the quantum affine algebra $\Uqsltwohat$.

\subsection{Some Hall numbers for $\mathcal A$}\label{ss:HallNumCohShP1}
We start with the following easy combinatorial lemma.
\begin{lemma}\label{le:FuncPhi}
The number $\varphi(a,b)$ of pairs $(J,L)\in\Fq[T,U]$ consisting of coprime homogeneous
polynomials of degree $a$ and $b$ respectively, is given by
\begin{equation*}
\varphi(a,b)=\begin{cases}(q-1)(q^{a+b+1}-1)&\text{if $a=0$ or $b=0$,}\\[2pt]
(q-1)(q^2-1)q^{a+b-1}&\text{if $a\geq1$ and $b\geq1$.}\end{cases}
\end{equation*}
\end{lemma}
\begin{proof}
Let $S$ be the set of pairs $(J,L)\in\Fq[T,U]$ consisting of non-zero homogeneous polynomials
of degree $a$ and $b$ respectively. The cardinality of $S$ is $(q^{a+1}-1)(q^{b+1}-1)$.
One can also count the number of elements in $S$ by factoring out a g.c.d.\ $D$ of $J$ and $L$.
For a fixed degree $d\leq\min(a,b)$, there are $(q^{d+1}-1)/(q-1)$ possibilities for $D$ up
to a non-zero scalar, and we thus get the relation
\begin{equation*}
(q^{a+1}-1)(q^{b+1}-1)=\sum_{d=0}^{\min(a,b)}{q^{d+1}-1\over q-1}\;\varphi(a-d,b-d).
\end{equation*}
The lemma then follows by induction on $\min(a,b)$.
\end{proof}

For any closed point $x$ of $\Pone(\Fq)$, let $q_x=q^{\deg x}$ be the cardinal of the
residue field of the local ring $\mathcal O_{\Pone(\Fq),x}$. We denote the greatest integer
less than or equal to a real number $a$ by~$\lfloor a\rfloor$. The following theorem
provides Hall numbers for the category $\mathcal A$.
\begin{theorem}\label{th:HallNumCohShP1}
In the Hall algebra $H(\mathcal A)$, one has the following relations:
\begin{enumerate}
\item\label{it:ThHNCSPa}$[\mathcal O(m)^{\oplus a}][\mathcal O(m)^{\oplus b}]=\left(\prod_{c=0}^{a-1}
{q^{a+b-c}-1\over q^{a-c}-1}\right)[\mathcal O(m)^{\oplus(a+b)}]$ for every
$m\in\mathbf Z$ and $a,b\in\mathbf N$.
\item\label{it:ThHNCSPb}If $\mathcal F=\mathcal O(n_1)\oplus\cdots\oplus\mathcal O(n_r)$ is a locally
free sheaf, if $m\in\mathbf Z$ is strictly greater than $n_1,\ldots,n_r$, and if $a$
is a non-negative integer, then
$[\mathcal F][\mathcal O(m)^{\oplus a}]=[\mathcal F\oplus\mathcal O(m)^{\oplus a}]$.
\item\label{it:ThHNCSPc}If $m<n$, then
\begin{equation*}
[\mathcal O(n)][\mathcal O(m)]=q^{n-m+1}[\mathcal O(m)\oplus\mathcal O(n)]+\sum_{a=1}^{
\lfloor(n-m)/2\rfloor}(q^2-1)\;q^{n-m-1}\;[\mathcal O(m+a)\oplus\mathcal O(n-a)].
\end{equation*}
\item\label{it:ThHNCSPd}If $\mathcal F$ and $\mathcal G$ are torsion sheaves whose support are disjoint,
then $[\mathcal F][\mathcal G]=[\mathcal F\oplus\mathcal G]$.
\item\label{it:ThHNCSPe}If $\mathcal F$ is a locally free sheaf and $\mathcal G$ is a torsion sheaf,
then $[\mathcal F][\mathcal G]=[\mathcal F\oplus\mathcal G]$.
\item\label{it:ThHNCSPf}If $x$ is a closed point, $r$ a positive integer and $n\in\mathbf Z$, then
\begin{equation*}
[\mathcal O_{(1^r)[x]}][\mathcal O(n)]=[\mathcal O(n+\deg x)\oplus
\mathcal O_{(1^{r-1})[x]}]+q_x^r\;[\mathcal O(n)\oplus\mathcal O_{(1^r)[x]}].
\end{equation*}
\end{enumerate}
\end{theorem}
\begin{proof}
\begin{enumerate}
\item Since the extension group $\Ext^1_{\mathcal A}(\mathcal O(m),\mathcal O(m))$ vanishes
by Corollary~\ref{co:VanishExt}, any short exact sequence of the form
\begin{equation*}
0\longrightarrow\mathcal O(m)^{\oplus b}\longrightarrow\mathcal F\longrightarrow
\mathcal O(m)^{\oplus a}\longrightarrow0
\end{equation*}
necessarily splits, and the product $[\mathcal O(m)^{\oplus a}][\mathcal O(m)^{\oplus b}]$
in $H(\mathcal A)$ is a scalar multiple of\break$[\mathcal O(m)^{\oplus(a+b)}]$. It remains to
compute the corresponding Hall number. Since $\End_{\mathcal A}(\mathcal O(m))\simeq\Fq$
by Lemma~\ref{le:HomLFCohShP1}~\ref{it:LeHLFCSPb}, this Hall number is equal to the number
of vector subspaces of dimension $b$ in a vector space of dimension $a+b$ over $\Fq$,
namely to $\left(\prod_{c=0}^{a-1}{q^{a+b-c}-1\over q^{a-c}-1}\right)$.
\item By Corollary~\ref{co:VanishExt}, the extension groups $\Ext^1_{\mathcal A}
(\mathcal O(n_i),\mathcal O(m))$ vanish. Thus any short exact sequence of the form
\begin{equation*}
0\longrightarrow\mathcal O(m)^{\oplus a}\,{\buildrel f\over\longrightarrow}\,
\mathcal G\longrightarrow\mathcal F\longrightarrow0
\end{equation*}
splits, and the product $[\mathcal F][\mathcal O(m)^{\oplus a}]$ in $H(\mathcal A)$ is
a scalar multiple of $[\mathcal F\oplus\mathcal O(m)^{\oplus a}]$. Let us put
$\mathcal G=\mathcal F\oplus\mathcal O(m)^{\oplus a}$ in the above short exact sequence,
and write $f=h\oplus i$ according to this decomposition. Then $h=0$, because all spaces
$\Hom_{\mathcal A}(\mathcal O(m),\mathcal O(n_i))$ vanish by Lemma~\ref{le:HomLFCohShP1}~\ref{it:LeHLFCSPc}.
It follows that $i$ is an automorphism. The number of suitable embeddings $f:\mathcal O(m)^{\oplus a}
\to\mathcal G$ is therefore equal to $\left|\Aut_{\mathcal A}(\mathcal O(m)^{\oplus a})\right|$, and
the Hall number we are looking for is equal to $1$.
\item By Proposition~\ref{pr:ExtCS1}, any short exact sequence of the form
\begin{equation*}
0\longrightarrow\mathcal O(m)\,{\buildrel f\over\longrightarrow}\,\mathcal F
\,{\buildrel g\over\longrightarrow}\,\mathcal O(n)\longrightarrow0
\end{equation*}
either splits, in which case $\mathcal F\simeq\mathcal O(m)\oplus\mathcal O(n)$, or there
exists $1\leq a\leq\lfloor(n-m)/2\rfloor$ such that $\mathcal F\simeq\mathcal O(m+a)\oplus
\mathcal O(n-a)$.

In the first case, we write $f=h\oplus i$ and $g=j\oplus\ell$, where $h\in\End_{\mathcal A}(\mathcal O(m))$,
$i,j\in\Hom_{\mathcal A}(\mathcal O(m),\mathcal O(n))$, and $\ell\in\End_{\mathcal A}(\mathcal O(n))$.
Since $\Hom_{\mathcal A}(\mathcal O(n),\mathcal O(m))=0$ by Lemma~\ref{le:HomLFCohShP1}~\ref{it:LeHLFCSPc},
the existence of a left inverse of $f$ requires that $h$ be an automorphism. Similarly, $\ell$ is an automorphism.
The map $i$ may be arbitrarily chosen and then the map $j$ should be equal to $-\ell\circ i\circ h^{-1}$.
Thus the set of suitable pairs $(f,g)$ is in one-to-one correspondence with $\Aut_{\mathcal A}
({\mathcal O(m)})\times\Aut_{\mathcal A}({\mathcal O(n)})\times\Hom_{\mathcal A}(\mathcal O(m),\mathcal O(n))$,
and the desired Hall number is $\left|\Hom_{\mathcal A}(\mathcal O(m),\mathcal O(n))\right|=q^{n-m+1}$.
This yields the term $q^{n-m+1}[\mathcal O(m)\oplus\mathcal O(n)]$ in the Hall product.

In the second case, the number of epimorphisms $g:\mathcal F\to\mathcal O(n)$ such that
$\ker g\simeq\mathcal O(m)$ is $(q-1)(q^2-1)q^{n-m-1}$ by Proposition~\ref{pr:ExtCS1} \vspace{2pt}
and Lemma~\ref{le:FuncPhi}. Since $\left|\Aut_{\mathcal A}(\mathcal O(n))\right|=q-1$,
the Hall number $\phi\,{}_{[\mathcal O(n)],[\mathcal O(m)]}^{[\mathcal F]}$ \vspace{2pt}is
$(q^2-1)q^{n-m-1}$, whence the term $(q^2-1)\;q^{n-m-1}\;[\mathcal O(m+a)\oplus\mathcal O(n-a)]$
in the Hall product.

\item\kern-\labelsep and\ \refstepcounter{enumi}\labelenumi\kern\labelsep They follow
from the vanishing of both $\Ext^1_{\mathcal A}(\mathcal F,\mathcal G)$ and $\Hom_{\mathcal A}
(\mathcal G,\mathcal F)$ (see\break Lemma~\ref{le:HomTCohShP1}~\ref{it:LeHTCSPa} and \ref{it:LeHTCSPd}).
\item This follows from Lemma~\ref{le:HomTCohShP1}~\ref{it:LeHTCSPd} and \ref{it:LeHTCSPe},
Lemma~\ref{le:HomMiCohShP1}, and Proposition~\ref{pr:ExtCS2} with the same reasoning as
for \ref{it:ThHNCSPc}.
\end{enumerate}
\end{proof}

\begin{other}{Application}
Let $x$ be a closed point of $\Pone(\Fq)$ and let $n\in\mathbf Z$. Using the relations in
Theorem~\ref{th:HallNumCohShP1} and the associativity of the Hall product, one obtains after
some calculation
\begin{multline*}
[\mathcal O_{[x]}][\mathcal O(n)^{\oplus2}]=q_x[\mathcal O(n)\oplus\mathcal O(n+\deg x)]\\
+q_x\left(1-\frac1q\right)\sum_{a=1}^{\lfloor(\deg x)/2\rfloor}[\mathcal O(n+a)\oplus
\mathcal O(n+\deg x-a)]+q_x[\mathcal O(n)^{\oplus2}\oplus\mathcal O_{[x]}].
\end{multline*}
In particular, if $\deg x\geq2$, then the Hall number $\phi_{[\mathcal O_{[x]}],[\mathcal
O(n)^{\oplus2}]}^{[\mathcal O(n+a)\oplus\mathcal O(n+\deg x-a)]}$ \vspace{2pt}is equal to
$q^{\deg x-1}(q-1)$ for each $1\leq a\leq\deg x-1$. By an analysis similar to those of
Propositions~\ref{pr:ExtCS1} and \ref{pr:ExtCS2}, one may deduce from this the following
fact, which we found not easy to prove directly:
\trivlist\item\itshape
If $P\in\Fq[T]$ is an irreducible polynomial of degree $d\geq2$ and if $1\leq a\leq d-1$,
then there are exactly $q^{d-1}(q-1)^2$ quadruples $(H,I,J,L)\in\Fq[T]^4$ consisting of polynomials
of degree $a$, $d-a$, $a-1$, $d-a-1$, respectively, such that $HI-JL=P$.
\endtrivlist
\end{other}

\subsection{The Hall subalgebras $H(\mathcal A_{\{x\}})$ and $H(\mathcal A_\tor)$}\label{ss:HallAlgA_tor}
The information provided by Theorem~\ref{th:HallNumCohShP1} is not sufficient to compute all
products in $H(\mathcal A)$. For instance, the elements $[\mathcal O_{r[x]}]$ do
not appear in its statement. More generally, it remains to understand how one can express the
elements $[\mathcal O_{\lambda[x]}]$ in terms of the elements $[\mathcal O_{(1^r)[x]}]$.

Let us fix a closed point $x$ of $\Pone(\Fq)$. The subcategory $\mathcal A_{\{x\}}$ of $\mathcal A$
defined in Section~\ref{ss:TorsionSheaves} is $k$-linear, abelian, and satisfies Conditions (H1)--(H4) of
Section~\ref{ss:InitData}. We can therefore consider the Hall algebra $H(\mathcal A_{\{x\}})$: it is the
subalgebra of $H(\mathcal A)$ spanned by the isomorphism classes of objects in $\mathcal A_{\{x\}}$.
Note that there is no difference between the Hall product $\cdot$ and the Ringel product $*$
on $H(\mathcal A_{\{x\}})$ since the Euler form vanishes on $K(\mathcal A_{\{x\}})$
by Proposition~\ref{pr:GrotGroup}~\ref{it:PrGGb}.

To simplify the notation, we will set $\widehat h_{r,x}=\sum_{|\lambda|=r}[\mathcal O_{\lambda[x]}]$,
where the sum runs over all partitions of weight $r$. The ring of symmetric polynomials over
the ground ring $\widetilde{\mathbf Z}=\mathbf Z[v,v^{-1}]/(v^2-q)$ in a countable infinite set of
indeterminates will be denoted by $\Lambda$. We will follow the notations of~\cite{Macdonald}
and denote the complete symmetric functions, the elementary symmetric functions, and the
Hall-Littlewood polynomials by $h_r\in\Lambda$, $e_r\in\Lambda$, and $P_\lambda(t)\in\Lambda[t]$,
respectively (see Sections~I.2 and III.2 of~[\textit{loc.~cit.}]). The next statement shows in
particular that the algebra $H(\mathcal A_{\{x\}})$ is commutative.

\begin{proposition}\label{pr:IsomHallPolSym}{\upshape(\cite{Kapranov}, Proposition~2.3.5)}
\begin{enumerate}
\item\label{it:PrIHPSa}There is a ring isomorphism $\Psi_x:H(\mathcal A_{\{x\}})\to
\Lambda$ that sends the elements $\widehat h_{r,x}$, $[\mathcal O_{(1^r)[x]}]$,
and $[\mathcal O_{\lambda[x]}]$ of $H(\mathcal A_{\{x\}})$, respectively, to the elements $h_r$,
$q_x^{-r(r-1)/2}e_r$ and $q_x^{-n(\lambda)}P_\lambda(q_x^{-1})$ of $\Lambda$, respectively,
for any integer $r\geq1$ and any partition $\lambda$.
\item\label{it:PrIHPSb}The $\widetilde{\mathbf Z}$-algebra $H(\mathcal A_{\{x\}})$ is
a polynomial algebra on the set $\{\widehat h_{r,x}\mid r\geq1\}$, as well as
on the set $\{[\mathcal O_{(1^r)[x]}]\mid r\geq1\}$. The family $([\mathcal O_{r[x]}])_{r\geq1}$
consists of algebraically independent elements and generates the $\mathbf Q[v]/(v^2-q)$-algebra
$H(\mathcal A_{\{x\}})\otimes_{\mathbf Z}\mathbf Q$.
\end{enumerate}
\end{proposition}
\begin{proof}
The isomorphism between the category $\mathcal A_{\{x\}}$ and the category of
$\mathcal O_{\Pone(\Fq),x}$-modules of finite length gives rise to an isomorphism
between their Hall algebras. Thus $H(\mathcal A_{\{x\}})$ is isomorphic to the Hall
algebra studied in Chapters~II and III of~\cite{Macdonald}. Assertion~\ref{it:PrIHPSa}
therefore follows from Paragraphs~III~(3.4), III.3~Example~1~(2), III.4~Example~1,
and III~(2.8) of [\textit{loc.~cit.}].

It is well-known that $\Lambda$ is the $\widetilde{\mathbf Z}$-algebra of polynomials
either in the complete symmetric functions $h_r$ or in the elementary symmetric functions
$e_r$, for $r\geq1$ (see Statements~I~(2.4) and I~(2.8) of~[\textit{loc.~cit.}]). This fact
implies the first assertion in Statement~\ref{it:PrIHPSb}. The second one follows from
Statement~III~(2.16) of~[\textit{loc.~cit.}] and its proof.
\end{proof}

We now define three generating functions in $H(\mathcal A_{\{x\}})[[s]]$ by
\begin{align*}
\widehat H_x(s)&=\,1+\sum_{r\geq1}\widehat h_{r,x}\;s^{r\deg x}
\,=\,\sum_{\beta\in\Iso(\mathcal A_{\{x\}})}\beta\;s^{\deg\beta},\displaybreak[0]\\
\widehat E_x(s)&=\,1+\sum_{r\geq1}(-1)^r\;q_x^{r(r-1)/2}\;[\mathcal O_{(1^r)[x]}]
\;s^{r\deg x},\displaybreak[0]\\
\widehat Q_x(s)&=\,1+\sum_{r\geq1}\;(1-q_x^{-1})\;v^{r\deg x}\;
[\mathcal O_{r[x]}]\;s^{r\deg x}.
\end{align*}
\begin{lemma}\label{le:RelaSerGen}
\begin{enumerate}
\item\label{it:LeRSGa}The following relations hold in $H(\mathcal A_{\{x\}})[[s]]$:
\begin{equation*}
\widehat H_x(s)\;\widehat E_x(s)=1\quad\text{and}\quad\widehat Q_x(s)=\frac{\widehat H_x(sv)}
{\widehat H_x(s/v)}.
\end{equation*}
\item\label{it:LeRSGb}In $H(\mathcal A_{\{x\}})[[s]]$, one has
\begin{equation*}
\widehat Q_x(s)=\sum\limits_{r\geq0}\;\left|\Aut_{\mathcal A}(\mathcal O_{r[x]})\right|
\;v^{-r\deg x}\;[\mathcal O_{r[x]}]\;s^{r\deg x}.
\end{equation*}
\end{enumerate}
\end{lemma}
\begin{proof}
Following Paragraphs~I~(2.2) and I~(2.5) of \cite{Macdonald}, we define generating series
in $\Lambda[[s]]$ by
\begin{equation*}
H(s)=1+\sum_{r\geq1}h_r\;s^r\qquad\text{and}\qquad E(s)=1+\sum_{r\geq1}e_r\;s^r.
\end{equation*}
By Formulae~I~(2.6) and III~(2.10) in~[\textit{loc.~cit.}], we have in $\Lambda[[s]]$
\begin{equation*}
H(s^{\deg x})E(-s^{\deg x})=1\quad\text{and}\quad 1+\sum_{r\geq1}(1-q_x^{-1})\;(sv)^{r\deg x}
\ P_{(r)}(q_x^{-1})=\frac{H((sv)^{\deg x})}{H((sv)^{\deg x}q_x^{-1})}.
\end{equation*}
Taking the inverse images by $\Psi_x$, one obtains the relations in Assertion~\ref{it:LeRSGa}.
As for Assertion~\ref{it:LeRSGb}, it follows from the equality $\left|\Aut_{\mathcal A}(\mathcal
O_{r[x]})\right|=q_x^r(1-q_x^{-1})$ (use Formula~II~(1.6) in~[\textit{loc.~cit.}] and
Proposition~\ref{pr:IsomCatA_x}).
\end{proof}

\begin{other}{Remark}
The category $\mathcal A_{\{x\}}$ satisfying Condition (H5), the algebra
$H(\mathcal A_{\{x\}})$ has the structure of a twisted Ringel-Green bialgebra. Since the
Euler form on $K(\mathcal A_{\{x\}})$ vanishes, the twist in the multiplication law on
$H(\mathcal A_{\{x\}})\otimes_{\widetilde{\mathbf Z}}H(\mathcal A_{\{x\}})$ is trivial,
so that $H(\mathcal A_{\{x\}})$ is a $\widetilde{\mathbf Z}$-bialgebra in the usual sense.
(This fact is due to Zelevinsky, see \cite{Green1}, p.~362; moreover, $H(\mathcal A_{\{x\}})$
has an antipode.) Now, $\Lambda$ is also a $\widetilde{\mathbf Z}$-bialgebra (see Example~25
in Section~I.5 of~\cite{Macdonald}). We claim that the isomorphism $\Psi_x$ defined
in Proposition~\ref{pr:IsomHallPolSym} preserves the coalgebra structures. To prove this, it
suffices to compare the behaviour of the coproduct of $H(\mathcal A_{\{x\}})$ on the
generators $\widehat h_{r,x}$ with that of the coproduct of $\Lambda$ on their images
$\Psi_x(\widehat h_{r,x})=h_r$. Using the definition of the coproduct in
Section~\ref{ss:HallAlgebra} and Proposition~\ref{pr:PropHallNum}~\ref{it:PrPHNf}, we
perform the following computation in~$(H(\mathcal A_{\{x\}})\otimes_{\widetilde{\mathbf Z}}
H(\mathcal A_{\{x\}}))[[s]]$:

\begin{align}
\Delta(\widehat H_x(s))&=\sum_{\beta\in\Iso(\mathcal A_{\{x\}})}s^{\deg\beta}\;\Delta(\beta)\nonumber\\
&=\sum_{\alpha,\beta,\gamma\in\Iso(\mathcal A_{\{x\}})}s^{\deg\beta}\;\frac{g_\alpha g_\gamma}
{g_\beta}\;\phi_{\alpha\gamma}^\beta\ (\alpha\otimes\gamma)\nonumber\\
&=\sum_{\alpha,\gamma\in\Iso(\mathcal A_{\{x\}})}\left(\sum_\beta\;\frac{g_\alpha g_\gamma}{g_\beta}
\;\phi_{\alpha\gamma}^\beta\right)\ (s^{\deg\alpha}\alpha)\otimes(s^{\deg\gamma}\gamma)\nonumber\\
&=\sum_{\alpha,\gamma\in\Iso(\mathcal A_{\{x\}})}(s^{\deg\alpha}\alpha)\otimes(s^{\deg\gamma}
\gamma)\nonumber\\&=\widehat H_x(s)\otimes\widehat H_x(s).\label{eq:CoprodHA_x}
\end{align}
Therefore $\Delta(\widehat h_{r,x})=\sum_{s=0}^r\widehat h_{s,x}\otimes\widehat h_{r-s,x}$ in
$H(\mathcal A_{\{x\}})\otimes_{\widetilde{\mathbf Z}}H(\mathcal A_{\{x\}})$. A similar
formula holds for the images of the complete symmetric functions $h_r$ by the coproduct of $\Lambda$
(see~[\textit{loc.~cit.}]), and our claim follows.
\end{other}

We now turn to the subcategory $\mathcal A_\tor$ of $\mathcal A$ consisting of all
torsion sheaves. The Hall algebra $H(\mathcal A_\tor)$ may be viewed as the subspace of
$H(\mathcal A)$ spanned by the isomorphism classes of objects in $\mathcal A_\tor$.
Since the category $\mathcal A_\tor$ is the direct sum of the categories $\mathcal A_{\{x\}}$
(see Lemma~\ref{le:HomTCohShP1}~\ref{it:LeHTCSPc}), the Hall algebra $H(\mathcal A_\tor)$ is
canonically isomorphic to the tensor product over $\widetilde{\mathbf Z}$ of the Hall
algebras $H(\mathcal A_{\{x\}})$ (this is Proposition~2.3.5~(a) in~\cite{Kapranov}).
Finally, we note that the Hall product $\cdot$ and the Ringel product $*$ coincide on
$H(\mathcal A_\tor)$ because the Euler form vanishes on $K(\mathcal A_\tor)$ by
Proposition~\ref{pr:GrotGroup}~\ref{it:PrGGb}.

We define elements $\widehat h_r$, $\widehat e_r$ and $\widehat q_r$ of $H(\mathcal A_\tor)$
for $r\geq1$ by means of the generating functions
\begin{align}
\widehat H(s)&=1+\sum_{r\geq1}\widehat h_rs^r=\prod_{x\in\Pone(\Fq)}\widehat H_x(s),
\label{eq:DefHhat}\displaybreak[0]\\
\widehat E(s)&=1+\sum_{r\geq1}\widehat e_rs^r=\prod_{x\in\Pone(\Fq)}\widehat E_x(s),
\nonumber\displaybreak[0]\\
\widehat Q(s)&=1+\sum_{r\geq1}\widehat q_rs^r=\prod_{x\in\Pone(\Fq)}\widehat Q_x(s).\nonumber
\end{align}
These equalities are meant to hold in $H(\mathcal A_\tor)[[s]]$; in the right-hand
side of the above equations, the products are over all closed points of $\Pone(\Fq)$.
\begin{lemma}\label{le:RelaHAtor}
\begin{enumerate}
\item\label{it:LeRHAa}One has the relations
\begin{equation*}
\widehat H(s)\;\widehat E(s)=1\quad\text{and}\quad\widehat Q(s)=\frac{\widehat H(sv)}
{\widehat H(s/v)}
\end{equation*}
or, equivalently, for each $r\geq1$,
\begin{gather*}
\widehat h_r+\sum_{s=1}^{r-1}\widehat h_s\;\widehat e_{r-s}+\widehat e_r=0,\\
(q^r-1)\;\widehat h_r=v^r\widehat q_r+\sum_{s=1}^{r-1}v^{r-s}\;\widehat h_s\;\widehat q_{r-s}.
\end{gather*}
\item\label{it:LeRHAb}The three families $(\widehat h_r)_{r\geq1}$, $(\widehat e_r)_{r\geq1}$, and
$(\widehat q_r)_{r\geq1}$ consist of algebraically independent elements.
\end{enumerate}
\end{lemma}
\begin{proof}
Assertion \ref{it:LeRHAa} follows from Lemma~\ref{le:RelaSerGen}~\ref{it:LeRSGa}.
Let $\Gamma$ be the subalgebra of $H(\mathcal A_\tor)$ generated by the subalgebras
$H(\mathcal A_{\{x\}})$ with $x\neq\infty$. Then $H(\mathcal A_\tor)$ is the algebra of
polynomials in the indeterminates $\widehat h_{r,\infty}$ with coefficients in $\Gamma$. It
is easy to see that $\widehat h_r-\widehat h_{r,\infty}$ belongs to $\Gamma[\widehat h_{1,\infty},
\cdots,\widehat h_{r-1,\infty}]$, which proves the algebraic independence of the elements $\widehat h_r$.
The algebraic independence of the other two families can then be deduced from Assertion~\ref{it:LeRHAa}.
This completes the proof of Assertion \ref{it:LeRHAb}.
\end{proof}

\subsection{A subalgebra of $H(\mathcal A)$}\label{ss:SubalgKapra}
The Ringel algebra $H(\mathcal A)$ turns out to be made of two parts: the first one is
the Ringel algebra $H(\mathcal A_\tor)$ described in Section~\ref{ss:HallAlgA_tor}, while
the second one is a certain subalgebra $B_1$ related to locally free sheaves. In this Section,
we explain this decomposition and use it to define a subalgebra $B$ of $H(\mathcal A)$
which will be related in Section~\ref{se:Uqsl2hat} to the quantum affine
algebra~$\Uqsltwohat$.

It will be necessary for us to extend the ground ring of the Ringel algebra $H(\mathcal A)$
and of certain $\widetilde{\mathbf Z}$-submodules $B$ of it from $\widetilde{\mathbf Z}$ to
a $\widetilde{\mathbf Z}$-algebra $R$. The $R$-module
$B\otimes_{\widetilde{\mathbf Z}}R$ will be denoted by $B_{(R)}$.

We first define the $q$-numbers, setting as usual $[a]=(v^a-v^{-a})/(v-v^{-1})$
for $a\in\mathbf Z$. We set $[a]!=\prod_{i=1}^a[i]$ for $a\geq1$, and agree that $[0]!=1$.
Remark that, up to a power of $v$, $[a]$ and $[a]!$ are non-zero integers.

We next record the following consequence of Theorem~\ref{th:HallNumCohShP1}
and of Proposition~\ref{pr:GrotGroup}.
\begin{lemma}\label{le:RelaHallLF}
\begin{enumerate}
\item\label{it:LeRHLFa}For all $m,n\in\mathbf Z$, one has
\begin{multline}
[\mathcal O(m+1)]*[\mathcal O(n)]-v^2\;[\mathcal O(n)]*[\mathcal O(m+1)]=\\
v^2\;[\mathcal O(m)]*[\mathcal O(n+1)]-[\mathcal O(n+1)]*[\mathcal O(m)].\label{eq:SubalgKEq1}
\end{multline}
\item\label{it:LeRHLFb}If $n_1<\cdots<n_r$ is an increasing sequence of integers and if
$c_1,\ldots,c_r$ is a sequence of positive integers, then one has
\begin{equation*}
[\mathcal O(n_1)]^{*c_1}*\cdots*[\mathcal O(n_r)]^{*c_r}=\left(\prod_{i=1}^r
q^{c_i(c_i-1)/2}[c_i]!\right)v^{\sum_{1\leq i<j\leq r}(n_j-n_i+1)c_ic_j}
\left[\bigoplus_{i=1}^r\mathcal O(n_i)^{\oplus c_i}\right].
\end{equation*}
\end{enumerate}
\end{lemma}
\begin{proof}
Assertion~\ref{it:LeRHLFa} can be proved by a tedious case by case examination,
using Relations~\ref{it:ThHNCSPa}--\ref{it:ThHNCSPc} in Theorem~\ref{th:HallNumCohShP1} and
Statement~\ref{it:PrGGb} of Proposition~\ref{pr:GrotGroup}. To prove Assertion~\ref{it:LeRHLFb},
one first compute
\begin{eqnarray*}
[\mathcal O(n_i)]^{*c_i}&=&v^{c_i(c_i-1)/2}\;[\mathcal O(n_i)]^{c_i}\\
&=&v^{c_i(c_i-1)/2}\;\left(\prod_{a=1}^{c_i}\frac{q^a-1}{q-1}\right)\;
[\mathcal O(n_i)^{\oplus c_i}]\\&=&q^{c_i(c_i-1)/2}\;[c_i]!\ [\mathcal
O(n_i)^{\oplus c_i}],
\end{eqnarray*}
using Theorem~\ref{th:HallNumCohShP1}~\ref{it:ThHNCSPa}, and then
\begin{align*}
[\mathcal O(n_1)]^{*c_1}*\cdots*\,&[\mathcal O(n_r)]^{*c_r}\\
&=\,v^{\sum_{1\leq i<j\leq r}(n_j-n_i+1)c_ic_j}\;
[\mathcal O(n_1)]^{*c_1}\cdots[\mathcal O(n_r)]^{*c_r}\\
&=\,\left(\prod_{i=1}^rq^{c_i(c_i-1)/2}[c_i]!\right)\;
v^{\sum_{1\leq i<j\leq r}(n_j-n_i+1)c_ic_j}\;
[\mathcal O(n_1)^{\oplus c_1}]\cdots[\mathcal O(n_r)^{\oplus c_r}]\\
&=\,\left(\prod_{i=1}^rq^{c_i(c_i-1)/2}[c_i]!\right)\;
v^{\sum_{1\leq i<j\leq r}(n_j-n_i+1)c_ic_j}\,
\left[\bigoplus_{i=1}^r\mathcal O(n_i)^{\oplus c_i}\right],
\end{align*}
using Theorem~\ref{th:HallNumCohShP1}~\ref{it:ThHNCSPb}.
\end{proof}

We now need two other pieces of notation. Let $C$ be the set of all sequences of
non-negative integers $\underline c=(c_n)_{n\in\mathbf Z}$ that have only finitely many
non-zero terms, and set
\begin{equation*}
X_{\underline c}=\prod_{n\in\mathbf Z}[\mathcal O(n)]^{*c_n}
\end{equation*}
for each $\underline c\in C$, the products being computed using the multiplication $*$
and the ascending order on $\mathbf Z$. We also denote the $\widetilde{\mathbf Z}$-submodule
of $H(\mathcal A)$ spanned by the isomorphism classes of locally free sheaves by $B_1$.

\begin{proposition}\label{pr:SubalgLF}
\begin{enumerate}
\item\label{it:PrSLFa}$B_1$ is a subalgebra of $H(\mathcal A)$.
\item\label{it:PrSLFb}If $R$ is a $\widetilde{\mathbf Z}$-algebra containing $\mathbf Q$,
then the family $(X_{\underline c})_{c\in C}$ is a basis of the $R$-module $(B_1)_{(R)}$.
\item\label{it:PrSLFc}The multiplication in the Ringel algebra $H(\mathcal A)$ induces an
isomorphism of $\widetilde{\mathbf Z}$-modules between $B_1\otimes_{\widetilde{\mathbf Z}}
H(\mathcal A_\tor)$ and $H(\mathcal A)$.
\end{enumerate}
\end{proposition}
\begin{proof}
Assertion \ref{it:PrSLFa} comes from the fact that $\mathcal A_\lf$ is a subcategory of
$\mathcal A$ closed under extensions (Proposition~\ref{pr:IndecCatA}~\ref{it:PrICAa})
and from the definition of the product $*$ in $H(\mathcal A)$. Assertion \ref{it:PrSLFb} follows from
Lemma~\ref{le:RelaHallLF}~\ref{it:LeRHLFb}. Finally, Proposition~\ref{pr:IndecCatA}~\ref{it:PrICAb}
and Theorem~\ref{th:HallNumCohShP1}~\ref{it:ThHNCSPe} imply Assertion~\ref{it:PrSLFc}.
\end{proof}

\begin{other}{Remark}\label{re:SmashProdStruct}
As mentioned in Section~\ref{ss:HallAlgA_tor}, the subalgebra $H(\mathcal A_\tor)$ of
$H(\mathcal A)$ is a Hopf algebra. Let us adopt Sweedler's notation and denote the image
of an element $a\in H(\mathcal A_\tor)$ under the coproduct by $\Delta(a)=\sum_{(a)}a_{(1)}
\otimes a_{(2)}$. The Hopf algebra $H(\mathcal A_\tor)$ acts on $H(\mathcal A)$ through
the adjoint representation, which is the homomorphism
$\ad:H(\mathcal A_\tor)\to\End_{\widetilde{\mathbf Z}}(H(\mathcal A))$ defined by
\begin{equation*}
a*x=\sum_{(a)}\left(\ad(a_{(1)})\cdot x\right)*a_{(2)},
\end{equation*}
for all $a\in H(\mathcal A_\tor)$ and $x\in H(\mathcal A)$. Relating the adjoint
action to Hecke operators, Kapranov has shown that $B_1$ is a $H(\mathcal A_\tor)$-submodule
of $H(\mathcal A)$ (see Proposition~4.1.1 in~\cite{Kapranov}). Assertion~\ref{it:PrSLFc} of
Proposition~\ref{pr:SubalgLF} can then be interpreted, in the language of Hopf algebras, as
stating that $H(\mathcal A)$ is the smash product of the Hopf algebra $H(\mathcal A_\tor)$
by the $H(\mathcal A_\tor)$-module algebra $B_1$.
\end{other}

The following result gives a commutation relation between certain elements of $B_1$ and
certain elements of $H(\mathcal A_\tor)$.
\begin{lemma}\label{le:RelaHallMixed}
For $n\in\mathbf Z$ and $r\geq1$, one has
\begin{equation}
\widehat h_r*[\mathcal O(n)]=\sum_{s=0}^r\;[s+1]\ [\mathcal O(n+s)]*\widehat h_{r-s}.
\label{eq:SubalgKEq2}
\end{equation}
\end{lemma}
\begin{proof}
We will use the generating series $\widehat E_x(s)$ and $\widehat H_x(s)$. We set
\begin{equation*}
X(t)=\sum_{n\in\mathbf Z}\;[\mathcal O(n)]\;t^n
\end{equation*}
and compute, using Relations \ref{it:ThHNCSPe} and \ref{it:ThHNCSPf} in Theorem
\ref{th:HallNumCohShP1}:
\begin{align*}
\widehat E_x(s)*X(t)&=\sum_{n\in\mathbf Z,\;r\geq0}\;(-1)^r\;s^{r\deg x}\;t^n\;q_x^{r(r-1)/2}
\;[\mathcal O_{(1^r)[x]}]*[\mathcal O(n)]\\
&=\sum_{n\in\mathbf Z,\;r\geq0}(-1)^r\;s^{r\deg x}\;t^n\;q_x^{r(r-1)/2}\\[-10pt]
&\hskip97pt\strut\smash{\times\left([\mathcal O(n)]*[\mathcal O_{(1^r)[x]}]+v^{(1-2r)\deg x}\;
[\mathcal O(n+\deg x)]*[\mathcal O_{(1^{r-1})[x]}]\right)}\\[8pt]
&=X(t)*\widehat E_x(s)\;\Bigl(1-(s/tv)^{\deg x}\Bigr).
\end{align*}
In view of Lemma~\ref{le:RelaSerGen}~\ref{it:LeRSGa}, we therefore have
\begin{equation*}
\widehat H(s)*X(t)=X(t)*\widehat H(s)\;\prod_{x\in\Pone(\Fq)}\frac1{1-(s/tv)^{\deg x}},
\end{equation*}
after expansion of the rational functions $1/(1-(s/tv)^{\deg x})$ in powers of $s/tv$.

Now in the formal power series ring $\mathbf Z[[s]]$, one has
\begin{equation*}
\prod_{x\in\Pone(\Fq)}{1\over1-s^{\deg x}}={1\over(1-s)(1-qs)},
\end{equation*}
where the product in the left hand side runs over all closed points of $\Pone(\Fq)$.
The previous equality follows from the calculation of the zeta function of $\Pone(\Fq)$
(see Section~C.1 of~\cite{Hartshorne} for a proof).
Therefore,
\begin{equation}
\widehat H(s)*X(t)=X(t)*\widehat H(s)\;\frac1{(1-s/tv)(1-sv/t)},\label{eq:RelaVertex}
\end{equation}
which is equivalent to our assertion.
\end{proof}

\begin{other}{Remark}
Lemma~\ref{le:RelaSerGen}~\ref{it:LeRSGb} shows that the elements $\psi_r\in H(\mathcal A_\tor)$
defined in Formula~(5.2) of~\cite{Kapranov} satisfy
\begin{equation*}
\widehat Q(s)=1+\sum_{r\geq1}v^{-r}\,\psi_r\,s^r.
\end{equation*}
On the other hand, using Relation (\ref{eq:RelaVertex}) above and
Lemma~\ref{le:RelaHAtor}~\ref{it:LeRHAa}, one obtains
\begin{equation*}
\widehat Q(s)*X(t)=X(t)*\widehat Q(s)\;\frac{(1-s/tq)}{(1-qs/t)}.
\end{equation*}
We thus recover Formula~(5.2.5) in~\cite{Kapranov}.
\end{other}

Finally, let $B_0$ be the subalgebra of $H(\mathcal A_\tor)$ generated by the family
$(\widehat h_r)_{r\geq1}$, and let $B$ be the subalgebra of the Ringel algebra
$H(\mathcal A_\tor)$ generated by $B_0$ and $B_1$. Let also $D$ be the set of all
sequences of non-negative integers $\underline d=(d_r)_{r\geq1}$ that have only
finitely many non-zero terms, and set
\begin{equation*}
\widehat h_{\underline d}=\prod_{r\geq1}\widehat h_r^{d_r},\quad
\widehat e_{\underline d}=\prod_{r\geq1}\widehat e_r^{d_r},\quad\text{and}\quad
\widehat q_{\underline d}=\prod_{r\geq1}\widehat q_r^{d_r},
\end{equation*}
for each $d\in D$.

\begin{proposition}\label{pr:SubalgKapra}
Let $R$ be a field of characteristic $0$ which is also a $\widetilde{\mathbf Z}$-algebra.
\begin{enumerate}
\item\label{it:PrSKa}The algebra $B_{(R)}$ is generated by the elements $[\mathcal O(n)]$
for $n\in\mathbf Z$ and the elements $\widehat h_r$ for~$r\geq1$.
\item\label{it:PrSKb}The families $(X_{\underline c}*\widehat h_{\underline d})_{(\underline c,
\underline d)\in C\times D}$, $(X_{\underline c}*\widehat e_{\underline d})_{(\underline c,
\underline d)\in C\times D}$, and $(X_{\underline c}*\widehat q_{\underline d})_{(\underline c,
\underline d)\in C\times D}$ are three bases of the $R$-module $B_{(R)}$.
\end{enumerate}
\end{proposition}
\begin{proof}
Proposition~\ref{pr:SubalgLF}~\ref{it:PrSLFc} implies that the multiplication $*$ induces
an isomorphism of $\widetilde{\mathbf Z}$-modules from $B_1\otimes_{\widetilde{\mathbf Z}}B_0$
to $B_1*B_0$. Lemma~\ref{le:RelaHallMixed} implies that $B_1*B_0$ is a subalgebra of the
Ringel algebra $H(\mathcal A)$, obviously equal to $B$. Since the
$\widetilde{\mathbf Z}$-modules $B_0$ and $B_1$ are free (see
Lemma~\ref{le:RelaHAtor}~\ref{it:LeRHAb}), the multiplication $*$ in the Ringel algebra
$H(\mathcal A)_{(R)}$ induces an isomorphism of $R$-modules from $(B_1)_{(R)}\otimes_R
(B_0)_{(R)}$ to $B_{(R)}$.

Lemma~\ref{le:RelaHAtor} shows that $(B_0)_{(R)}$ is a polynomial algebra on each of the
three set of indeterminates: $\{\widehat h_r\mid r\geq1\}$, $\{\widehat e_r\mid r\geq1\}$,
or $\{\widehat q_r\mid r\geq1\}$. Thus the families $(\widehat h_{\underline d})_{\underline
d\in D}$, $(\widehat e_{\underline d})_{\underline d\in D}$, and $(\widehat q_{\underline
d})_{\underline d\in D}$ are three bases of the $R$-vector space $(B_0)_{(R)}$.
Both assertions of our proposition follow now from Proposition~\ref{pr:SubalgLF}~\ref{it:PrSLFb}.
\end{proof}

\begin{other}{Remark}
In view of Remark~\ref{re:SmashProdStruct}, the fact that $B$ is a subalgebra of
$H(\mathcal A)$ should be considered as a consequence of the fact that $B_0$ is a
sub-bialgebra of $H(\mathcal A_\tor)$, itself a consequence of Formulae~(\ref{eq:CoprodHA_x})
and (\ref{eq:DefHhat}).
\end{other}

\section{Link with the quantum affine algebra $\Uqsltwohat$}\label{se:Uqsl2hat}
Our aim now is to describe the relationship between the Hall algebra $H(\mathcal A)$
investigated in Section~\ref{se:HallAlgCohShP1} and the quantum affine algebra
$\Uqsltwohat$. In the usual definitions of the latter, $q$ is an indeterminate.
It will however be more convenient for us to deal with a specialized version of
$\Uqsltwohat$, in which $q$ is the number of elements of the finite field that we have
chosen at the beginning of Section~\ref{se:HallAlgCohShP1}. We therefore fix for the remainder
of this paper a field $R$ of characteristic $0$ together with \vspace{2pt}a square root $v$
of the number $q$.

In this section, we first recall the definition of the $R$-algebra $\Uqsltwohat$ in its
loop-like realization and define a certain subalgebra $V^+$ in it. We then present an
elementary proof of Kapranov's result asserting that the $R$-algebra $V^+$ is isomorphic to
the $R$-algebra $B_{(R)}$ defined in Section~\ref{ss:SubalgKapra}. We end
with several comments, observing that Kapranov's approach to~$\Uqsltwohat$ sheds a new light
on certain recent constructions by Beck, Chari, and Pressley~\cite{Beck-Chari-Pressley}.

\subsection{Definition of $\Uqsltwohat$}\label{ss:DefUqsl2hat}
Following Drinfeld \cite{Drinfeld2}, we define $\Uqsltwohat$ as the $R$-algebra generated
by elements $K^{\pm1}$, $C^{\pm1/2}$, $h_r$, where $r\in\mathbf Z\setminus\{0\}$, and
$x_n^\pm$, where $n\in\mathbf Z$, submitted to the relations
\begin{xalignat}2
&K\;K^{-1}=K^{-1}\;K=1,\nonumber\displaybreak[0]\\[2pt]
&C^{1/2}\;C^{-1/2}=C^{-1/2}\;C^{1/2}=1,\nonumber\displaybreak[0]\\[2pt]
&C^{1/2}\quad\text{is central,}\nonumber\displaybreak[0]\\[2pt]
&[K,h_r]=0&&\text{for $r\in\mathbf Z\setminus\{0\}$,}\nonumber\displaybreak[0]\\[2pt]
&K\;x_n^\pm=v^{\pm2}x_n^\pm\;K&&\text{for $n\in\mathbf Z$,}\nonumber\displaybreak[0]\\[2pt]
&[h_r,h_s]=\delta_{r,-s}\;\frac{[2r]}r\;\frac{C^r-C^{-r}}{v-v^{-1}}&&\text{for $r,s\in\mathbf Z
\setminus\{0\}$,}\nonumber\displaybreak[0]\\[2pt]
&[h_r,x_n^\pm]=\pm\frac{[2r]}r\;C^{\mp|r|/2}x_{n+r}^\pm&&\text{for $n,r\in\mathbf Z$,
$r\neq0$,}\nonumber\displaybreak[0]\\[2pt]
&x_{m+1}^\pm\,x_n^\pm-v^{\pm2}x_n^\pm\,x_{m+1}^\pm=v^{\pm2}x_m^\pm\,x_{n+1}^\pm-
x_{n+1}^\pm\,x_m^\pm&&\text{for $m,n\in\mathbf Z$,}\label{eq:UqEq1}\displaybreak[0]\\[2pt]
&[x_m^+,x_n^-]=\frac{C^{(m-n)/2}\psi_{m+n}^+-C^{(n-m)/2}\psi_{m+n}^-}{v-v^{-1}}
&&\text{for $m,n\in\mathbf Z$,}\nonumber
\end{xalignat}
where the elements $\psi_{\pm r}^{\pm}$ are defined by the generating functions
\begin{equation*}
\sum_{r\geq0}\psi_{\pm r}^\pm\,s^{\pm r}=K\exp\left(\pm(v-v^{-1})\sum_{r\geq1}h_{\pm r}\,
s^{\pm r}\right)
\end{equation*}
for $r\geq0$ and are defined to be zero for $r\leq-1$.

Relying in part on previous work of Damiani~\cite{Damiani}, Beck \cite{Beck} made precise
the link between this definition and Drinfeld's and Jimbo's original definition
\cite{Drinfeld1,Jimbo} of $\Uqsltwohat$ as the quantized enveloping algebra
associated to the generalized Cartan matrix $\begin{pmatrix}2&-2\\-2&2\end{pmatrix}$ of
type $A_1^{(1)}$.

Let $V^+$ be the subalgebra of $\Uqsltwohat$ generated by the elements $x_n^+$ and $h_r\,C^{r/2}$,
for $n\in\mathbf Z$ and $r\geq1$. The aim of Section~\ref{se:Uqsl2hat} is to prove the following result.
\begin{theorem}\label{th:ThmKapranov}{\upshape(\cite{Kapranov}, Theorem~5.2.1)}
The $R$-algebras $B_{(R)}$ and $V^+$ are isomorphic.
\end{theorem}

\subsection{Structure of $\Uqsltwohat$}\label{ss:StructUq}
Following Section~1 of~\cite{Beck-Chari-Pressley}, we define elements
$\widetilde\psi_{\pm r}^\pm$ for $r\geq1$ by
\begin{equation}
1\pm\sum_{r\geq1}(v-v^{-1})\,\widetilde\psi_{\pm r}^\pm\,s^{\pm r}\;=\;
\exp\left(\pm(v-v^{-1})\sum_{r\geq1}h_{\pm r}\,C^{\pm r/2}\,s^{\pm r}\right).\label{eq:Defpsitilde}
\end{equation}

\pagebreak[3]
\noindent Let us denote by
\begin{itemize}
\item$N^\pm$ the subalgebra generated by the elements $x_n^\pm$, where $n\in\mathbf Z$;
\item$H$ the subalgebra generated by the elements $K^{\pm1}$, $C^{\pm1/2}$, and
$h_r$, where $r\in\mathbf Z\setminus\{0\}$;
\item$H^\pm$ the subalgebra generated by the elements $\widetilde\psi_{\pm r}^\pm$, where $r\geq1$;
\item$H^0$ the subalgebra generated by the elements $K^{\pm1}$ and $C^{\pm1/2}$.
\end{itemize}

\begin{proposition}\label{pr:DecTriangUq}\
\begin{enumerate}
\item\label{it:LeSUa}The multiplication induces a linear isomorphism
$N^-\otimes_RH\otimes_RN^+\to\Uqsltwohat$.
\item\label{it:LeSUb}The multiplication induces a linear isomorphism
$H^-\otimes_RH^0\otimes_RH^+\to H$.
\item\label{it:LeSUc}The generators $\widetilde\psi_r^+$ ($r\geq1$) of the algebra $H^+$ are
algebraically independent.
\item\label{it:LeSUd}The family of products $\left(\prod_{n\in\mathbf Z}(x_n^+)^{c_n}
\right)_{\underline c\in C}$, performed in the ascending order of $\mathbf Z$, is a basis of $N^+$.
\end{enumerate}
\end{proposition}
\begin{proof}
It is asserted in Proposition~12.2.2 of~\cite{Chari-Pressley} that the map $N^-\otimes_RH\otimes_RN^+
\to\Uqsltwohat$ induced by the multiplication of $\Uqsltwohat$ is surjective. The defining
relations of $\Uqsltwohat$ imply easily that the map $H^-\otimes_RH^0\otimes_RH^+\to H$ induced by
the multiplication of $\Uqsltwohat$ is also surjective. The algebra $H^\pm$ is generated by the
pairwise commuting elements $\widetilde\psi_{\pm r}^\pm$ for $r\geq1$, which shows that the
monomials $\left(\prod_{r\geq1}(\widetilde\psi_r^\pm)^{d_r}\right)$, for $\underline d\in D$,
span the $R$-vector space $H^\pm$. Similarly, the family of elements \vspace{2pt}
$(K^aC^{b/2})_{(a,b)\in\mathbf Z^2}$ span the $R$-vector space $H^0$. Finally, an easy induction
shows that the products $\left(\prod_{n\in\mathbf Z}(x_{\pm n}^\pm)^{c_n}\right)$, \vspace{2pt}
performed in the ascending order of $\mathbf Z$ and for $\underline c\in C$, span the $R$-vector
space $N^\pm$. Consequently, the elements
\begin{equation*}
M(a,b,\underline{c'},\underline{c''},\underline{d'},\underline{d''})=
\left(\prod_{n\in\mathbf Z}(x_{-n}^-)^{c'_n}\right)
\left(\prod_{r\geq1}(\widetilde\psi_{-r}^-)^{d'_r}\right)
K^aC^{b/2}
\left(\prod_{r\geq1}(\widetilde\psi_r^+)^{d''_r}\right)
\left(\prod_{n\in\mathbf Z}(x_n^+)^{c''_n}\right),
\end{equation*}
where $a,b\in\mathbf Z$, $\underline{c'},\underline{c''}\in C$, and $\underline{d'},
\underline{d''}\in D$, span the $R$-vector space $\Uqsltwohat$.

Now observe that the definition of $\Uqsltwohat$ implies the existence of an automorphism $T$ of
the $R$-algebra $\Uqsltwohat$ such that
\begin{equation*}
T(x_n^\pm)=x_{n\mp1}^\pm,\quad T(K^{\pm1})=K^{\pm1},\quad T(C^{\pm1/2})=C^{\pm1/2},\quad
T(h_r)=h_r,
\end{equation*}
for all $n,r\in\mathbf Z$ with $r\neq0$. (Using Proposition~3.10.2~(b) and Definition~4.6 of
\cite{Beck}, one can easily see that $T$ is the automorphism of $\Uqsltwohat$ that lifts the
translation along the fundamental weight to the braid group of the extended affine Weyl group
of $\mathfrak{sl}_2$.) On the other hand, observe, as a consequence of the
Poincar\'e-Birkhoff-Witt theorem for $\Uqsltwohat$ (Proposition~6.1 in~\cite{Beck}), that
the set of elements
\begin{equation*}
\bigl\{M(a,b,\underline{c'},\underline{c''},\underline{d'},\underline{d''})\bigm|
a,b\in\mathbf Z,\;\underline{c'},\underline{c''}\in C^2,\;\underline{d'},\underline{d''}\in D^2,
\;n<0\Rightarrow c'_n=c''_n=0\bigr\}
\end{equation*}
is linearly independent over $R$. Using this and the automorphism $T$,
one proves the linear independence over $R$ of the family of elements
$(M(a,b,\underline{c'},\underline{c''},\underline{d'},\underline{d''}))_{(a,b,\underline{c'},
\underline{c''},\underline{d'},\underline{d''})\in\mathbf Z^2\times C^2\times D^2}$.

This family is therefore a basis of $\Uqsltwohat$, which entails simultaneously all the
assertions of the lemma.
\end{proof}

Let us remark that the algebra $H^+$ is the subalgebra denoted by $\mathbf U^+(0)$
in~\cite{Beck-Chari-Pressley} (see Proposition~1.3~(iii) in~[\textit{loc.~cit.}] for instance).
Following Section~1 of [\textit{loc.~cit.}], we now define elements $\widetilde P_r$ and $P_r$ of $H^+$,
for $r\geq1$, by the following generating functions:
\begin{eqnarray}
\widetilde P(s)=1+\sum_{r\geq1}\widetilde P_rs^n&=&\exp\left(\sum_{r\geq1}\frac{h_r\,C^{r/2}}{[r]}
\;s^r\right),\label{eq:DefP_itilde}\\P(s)=1+\sum_{r\geq1}P_rs^r&=&\exp\left(-\sum_{r\geq1}
\frac{h_r\,C^{r/2}}{[r]}\;s^r\right).\label{eq:DefP_i}
\end{eqnarray}

For sequences $\underline c=(c_n)_{n\in\mathbf Z}\in C$ and $\underline d=(d_r)_{r\geq1}
\in D$, we define
\begin{equation*}
x_{\underline c}^+=\prod_{n\in\mathbf Z}x_n^{c_n},\quad
\widetilde P_{\underline d}=\prod_{r\geq1}\widetilde P_r^{d_r},\quad
P_{\underline d}=\prod_{r\geq1}P_r^{d_r},\quad\text{and}\quad
\widetilde\psi_{\underline d}^+=\prod_{r\geq1}(\widetilde\psi_r^+)^{d_r}.
\end{equation*}

\begin{proposition}\label{pr:SubalgUqPlus}
\begin{enumerate}
\item\label{it:PrUPa}The algebra $V^+$ is generated by the elements $x_n^+$ and $\widetilde P_r$,
for $n\in\mathbf Z$ and $r\geq1$.
\item\label{it:PrUPb}The families $(x_{\underline c}^+\widetilde P_{\underline d}^{})_{(\underline c,
\underline d)\in C\times D}^{}$, $(x_{\underline c}^+P_{\underline d}^{})_{(
\underline c,\underline d)\in C\times D}^{}$, and $(x_{\underline c}^+\widetilde
\psi_{\underline d}^+)_{(\underline c,\underline d)\in C\times D}^{}$ are three
bases of the $R$-vector space $V^+$.
\end{enumerate}
\end{proposition}
\begin{proof}
By definition, $V^+$ is the subalgebra of $\Uqsltwohat$ generated by the elements $x_n^+$ and
$h_rC^{r/2}$, for $n\in\mathbf Z$ and $r\geq1$. Assertion~\ref{it:PrUPa} follows therefore from
the definition of the elements $\widetilde P_r$ (formula~(\ref{eq:DefP_itilde})) and from
the fact that the scalars $[r]$ do not vanish in the field $R$ for any $r\geq1$. 
 
Proposition~2.8 in~\cite{Beck-Chari-Pressley} states that for all integers $n\geq0$ and
$r\geq 1$, one has
\begin{equation}
\widetilde P_r\;x_n^+=\sum_{s=0}^r\;[s+1]\;x_{n+s}^+\;\widetilde P_{r-s}.\label{eq:UqEq2}
\end{equation}
Applying a well-chosen power of the automorphism $T$ defined in the proof of
Proposition~\ref{pr:DecTriangUq}, one immediately sees that Formula~(\ref{eq:UqEq2})
holds more generally for each $n\in\mathbf Z$. Together with
Proposition~\ref{pr:DecTriangUq}~\ref{it:LeSUa}, this shows that the multiplication
map induces a linear isomorphism $N^+\otimes_RH^+\to V^+$.

By Formulae~(1.8) and (1.9) in~\cite{Beck-Chari-Pressley}, we have the following
relations in $H^+[[s]]$:
\begin{equation*}
\widetilde P(s)\;P(s)=1\quad\text{and}\quad\frac{\widetilde P(sv)}{\widetilde P(s/v)}=1+\sum_{r\geq1}
(v-v^{-1})\;\widetilde\psi_r^+\;s^r,
\end{equation*}
or, equivalently, for each $r\geq1$:
\begin{gather*}
\widetilde P_r+\sum_{s=1}^{r-1}\widetilde P_sP_{r-s}+P_r=0,\\
[r]\,\widetilde P_r=\widetilde\psi_r^++\sum_{s=1}^{r-1}v^{-s}\;\widetilde P_s\;\widetilde\psi_{r-s}^+.
\end{gather*}
Together with Assertion~\ref{it:LeSUc} of Proposition~\ref{pr:DecTriangUq}, this implies that $H^+$
is a polynomial algebra on each of the three set of indeterminates: $\{\widetilde\psi_r\mid r\geq1\}$,
$\{\widetilde P_r\mid r\geq1\}$, or $\{P_r\mid r\geq1\}$. Thus the families
$(\widetilde P_{\underline d})_{\underline d\in D}$, $(P_{\underline d})_{\underline d\in D}$,
and $(\widetilde\psi_{\underline d}^+)_{\underline d\in D}$ are three bases of the $R$-vector
space $H^+$. Assertion~\ref{it:PrUPb} follows from this and from
Proposition~\ref{pr:DecTriangUq}~\ref{it:LeSUd}.
\end{proof}

Theorem~\ref{th:ThmKapranov} is now evident. The isomorphism sends $[\mathcal O(n)]$
to $x_n^+$, $\widehat h_r$ to $\widetilde P_r$, $\widehat e_r$ to $P_r$, and
$\widehat q_r$ to $(v-v^{-1})\widetilde\psi_r^+$, respectively. Relations~(\ref{eq:SubalgKEq1})
and (\ref{eq:SubalgKEq2}) correspond to Relations (\ref{eq:UqEq1}) and (\ref{eq:UqEq2}). 

\subsection{Concluding remarks}\label{ss:ConclRems}
As mentioned in the introduction, Ringel was the first one to discover relations
between Hall algebras and quantized enveloping algebras. In \cite{Ringel4} he noticed that,
in his context, the natural basis of the Hall algebra corresponds to a basis of type
Poincar\'e-Birkhoff-Witt of Lusztig's integral form of the positive part of the quantized enveloping
algebra. Here a similar phenomenon occurs:
\begin{itemize}
\item by Lemma~\ref{le:RelaHallLF}~\ref{it:LeRHLFb}, the element $[\mathcal O(n_1)^{\oplus c_1}
\oplus\cdots\oplus\mathcal O(n_r)^{\oplus c_r}]$ of $H(\mathcal A)$ is equal, up to a power
of $v$, to the product of divided powers
\begin{equation*}
\left(\frac1{[c_1]!}[\mathcal O(n_1)]^{*c_1}\right)*\cdots*\left(\frac1{[c_r]!}
[\mathcal O(n_r)]^{*c_r}\right);
\end{equation*}
\item by Lemma~2.3 in~\cite{Beck-Chari-Pressley}, the monomials in the $\widehat h_r$ for $r\geq1$
correspond to the elements of a Poincar\'e-Birkhoff-Witt basis of Lusztig's integral form of $H^+$.
\end{itemize}
These observations are likely to be part of a more complete statement, for which one would need
a version of the Hall algebra $H(\mathcal A)$ (or at least of the algebra $B$) with a generic
parameter $q$ as well as an integral version of \vspace{2pt}Proposition~\ref{pr:DecTriangUq}.

The elements $P_r$, $\widetilde P_r$, and $\widetilde\psi_r$ play important r\^oles for $\Uqsltwohat$,
namely in the classification of the finite-dimensional simple $\Uqsltwohat$-modules (see Theorem~12.2.6
in~\cite{Chari-Pressley}) and in the construction of a Poincar\'e-Birkhoff-Witt basis of Lusztig's
integral form of $\Uqsltwohat$ (see Theorem~2 of~\cite{Beck-Chari-Pressley}). The relations
(\ref{eq:Defpsitilde}), (\ref{eq:DefP_itilde}), and (\ref{eq:DefP_i}) defining them, though explicit,
look rather artificial. Kapranov's isomorphism $B_{(R)}\to V^+$ helps show where these elements come from:
they are the standard generators of the algebra of symmetric polynomials (elementary symmetric functions,
complete symmetric functions, Hall-Littlewood polynomials) carried over to the Hall algebras
$H(\mathcal A_{\{x\}})$ and averaged over the closed points of~$\Pone(\Fq)$.

\bigskip
\noindent Institut de Recherche Math\'ematique Avanc\'ee\\
Universit\'e Louis Pasteur et CNRS\\
7, rue Ren\'e Descartes\\
F-67084 Strasbourg Cedex\\
France\\[6pt]
E-mail: \texttt{baumann@math.u-strasbg.fr}, \texttt{kassel@math.u-strasbg.fr}\\
Fax: (+33/0) 3 88 61 90 69\\
Web: \texttt{http://www-irma.u-strasbg.fr/\~{}kassel/}
\end{document}